\documentclass[onespacing]{elsart}
\usepackage{amstext,amssymb,amsmath,amscd}
\usepackage{graphics,graphicx,epsfig,subfigure,enumerate,rotating,theorem,algorithm,algorithmic,array,longtable,multirow,multicol}
\usepackage[numbers]{natbib}
\numberwithin{equation}{section}
\theoremstyle{plain}
\newtheorem{theorem}{Theorem}[section]
\newtheorem{lemma}{Lemma}[section]

\newtheorem{remark}{Remark}[section]
\newtheorem{proposition}{Proposition}[section]

\newcommand{\psnr}{\mathrm{PSNR}}
\newcommand{\psnrin}{$\psnr_{\mathrm{in}}$}
\newcommand{\card}{\mathrm{Card}}
\newcommand{\indic}[1]{\mathbf{1}_{#1}}

\begin{document}

\begin{frontmatter}



\title{Stein Block Thresholding For Image Denoising}

\author[1]{C. Chesneau\corauthref{cor1}}
\author[2]{J. Fadili}
\author[3]{J.-L. Starck}
\address[1]{Laboratoire de Math\'ematiques Nicolas Oresme, CNRS-Universit\'{e} de Caen, Campus II, Science 3, 14032, Caen Cedex, France.}
\address[2]{GREYC CNRS-ENSICAEN-Universit\'{e} de Caen, Image Processing Group, 14050, Caen Cedex, France.}
\address[3]{Laboratoire AIM, CEA/DSM-CNRS-Universit\'{e} Paris Diderot, IRFU, SEDI-SAP, Service d'Astrophysique, Centre de Saclay, 91191 Gif-Sur-Yvette cedex, France.}
\corauth[cor1]{Corresponding author: chesneau@math.unicaen.fr}

\begin{abstract}
In this paper, we investigate the minimax properties of Stein block thresholding in any dimension $d$ with a particular emphasis on $d=2$. Towards this goal, we consider a frame coefficient space over which minimaxity is proved. The choice of this space is inspired by the characterization provided in \cite{BorupNielsen} of family of smoothness spaces on $\mathbb{R}^d$, a subclass of so-called decomposition spaces \cite{Feichtinger}. These smoothness spaces cover the classical case of Besov spaces, as well as smoothness spaces corresponding to curvelet-type constructions. Our main theoretical result investigates the minimax rates over these decomposition spaces, and shows that our block estimator can achieve the optimal minimax rate, or is at least nearly-minimax (up to a $\log$ factor) in the least favorable situation. Another contribution is that the minimax rates given here are stated for a general noise sequence model in the transform coefficient domain beyond the usual i.i.d. Gaussian case. The choice of the threshold parameter is theoretically discussed and its optimal value is stated for some noise models such as the (non-necessarily i.i.d.) Gaussian case. We provide a simple, fast and a practical procedure. We also report a comprehensive simulation study to support our theoretical findings. The practical performance of our Stein block denoising compares very favorably to the BLS-GSM state-of-the art denoising algorithm on a large set of test images. A toolbox is made available for download on the Internet to reproduce the results discussed in this paper.
\end{abstract}


\end{frontmatter}

\section{Introduction}

Consider the nonparametric regression model:
\begin{equation}\label{mod}
Y_{\bf i}=f({\bf i}/n)+ \epsilon_{\bf i}, \ \ \ \ \ \ {\bf i}\in \{1, ..., n\}^d,
\end{equation}
where $d\in \mathbb{N}^*$ is the dimension of the data, $(Y_{\bf i})_{{\bf i}\in \{1, ..., n\}^d}$ are the observations regularly sampled on a $d$-dimensional Cartesian grid, $(\epsilon_{\bf i})_{{\bf i}\in \{1, ..., n\}^d}$ are independent and identically distributed (i.i.d.) $\mathcal{N}(0,1)$, and $f : [0,1]^d \rightarrow \mathbb{R}$ is an unknown function. The goal is to estimate $f$ from the observations.
We want to build an adaptive estimator $\widehat f$ (i.e. its construction depends on the observations only) such that the mean integrated squared error (MISE) defined by
$R(\widehat f,f)=\mathbb{E}\left( \int_{[0,1]^d}\left( \widehat f({\bf x}) -f({\bf x})\right)^2d{\bf x}\right)$ is as small as possible for a wide class of $f$.
A now classical approach to the study of nonparametric problems of the form \eqref{mod} is to, first, transform the data to obtain a sequence of coefficients, second, analyze and process the coefficients (e.g. shrinkage, thresholding), and finally, reconstruct the estimate from the processed coefficients. This approach has already proven to be very successful by several authors and a good survey may be found in \cite{TsybakovBook,Johnstone99,JohnstoneMonograph}. In particular, it is now well established that the quality of the estimation is closely linked to the sparsity of the sequence of coefficients representing $f$ in the transform domain. Therefore, in this paper, we focus our attention on transform-domain shrinkage methods, such as those operating in the wavelet domain.

\subsection{The one-dimensional case}
First of all, let's consider the one-dimensional case $d=1$. The most standard of wavelet shrinkage methods is VisuShrink of \citep{donohoj}. It is constructed through individual (or term-by-term) thresholding
of the empirical wavelet coefficients. It enjoys good theoretical (and practical) properties. In particular, it achieves the optimal rate of convergence up to a logarithmic term over the H\"older class under the MISE.
In other words, if $\widehat f^{V} $ denotes VisuShrink, and $\Lambda^s(M)$ the H\"older smoothness class, then there exists a constant $C>0$ such that
\begin{equation}\label{eq1}
\sup_{f\in \Lambda^s(M)}R(\widehat f^{V}, f)\le C n^{-2s/(1+2s)} (\log n)^{2s/(1+2s)} .
\end{equation}
Other term-by-term shrinkage rules have been developed. See, for instance, the firm shrinkage \citep{gaobruce} or the non-negative garrote shrinkage \citep{gao}. In particular, they satisfy \eqref{eq1} but improve the value of the constant $C$. An exhaustive account of other shrinkage methods is provided in \citep{antoniadis} that the interested reader may refer to.

The individual approach achieves a degree of trade-off between variance and bias contribution to the MISE. However, this trade-off is not optimal; it removes too many terms from the observed wavelet expansion, with the consequence the estimator is too biased and has a sub-optimal MISE convergence rate (and also in other $L_p$ metrics $1 \leq p \leq \infty$). One way to increase estimation precision is by exploiting information about neighboring coefficients. In other words, empirical wavelet coefficients tend to form clusters that could be thresholded in blocks (or groups) rather than individually. This would allow threshold decisions to be made more accurately and permit convergence rates to be improved. Such a procedure has been introduced in \citep{hall1,hall2} who studied wavelet shrinkage methods based on block thresholding.
The procedure first divides the wavelet coefficients at each resolution
level into non-overlapping blocks and then keeps all
the coefficients within a block if, and only if, the magnitude of the sum of the squared
empirical coefficients within that block is greater than a fixed threshold.
The original procedure developed by \citep{hall1,hall2} is defined with the block size $(\log n)^2$.
BlockShrink of \citep{cai97,cai02} is the optimal version of this procedure. It uses a different block size, $\log n$, and enjoys a number of advantages over the conventional individual
thresholding. In particular, it achieves the optimal rate of convergence over the H\"older class under the MISE.
In other words, if $\widehat f^{B}$ denotes the BlockShrink estimate, then there exists a constant $C>0$ such that
\begin{equation}\label{eq2}
\sup_{f\in \Lambda^s(M)}R(\widehat f^{B}, f)\le C n^{-2s/(1+2s)} .
\end{equation}
Clearly, in comparison to VisuShrink, BlockShrink removes the extra logarithmic term. The minimax properties of BlockShrink under the $L_p$ risk have been studied in \citep{Chesneau08b}. Other local block thresholding rules have been developed. Among them, there is BlockJS of \citep{cai99,cai02} which combines James-Stein rule (see \citep{stein}) with the wavelet methodology.
In particular, it satisfies \eqref{eq2} but improves the value of the constant $C$. From a practical point view, it is better than BlockShrink. Further details about the theoretical performances of BlockJS can be found in \citep{cavalier}. We refer to \citep{antoniadis} and \citep{caisil} for a comprehensive simulation study. Variations of BlockJS are BlockSure of \citep{chicken} and SureBlock of \citep{cai07}. The distinctive aspect of these block thresholding procedures is to provide data-driven algorithms to chose the threshold parameter. Let's also mention the work of \cite{Abramovich02} who considered wavelet block denoising in a Bayesian framework to obtain level-dependent block shrinkage and thresholding estimates.

\subsection{The multi-dimensional case}
Denoising is a long-standing problem in image processing. Since the seminal papers by Donoho \& Johnstone \citep{donohoj}, the image processing literature has been inundated by hundreds of papers applying or proposing modifications of the original algorithm in image  denoising. Owing to recent advances in computational harmonic analysis, many multi-scale geometrical transforms, such as ridgelets \citep{CandesDonohoRidgelets}, curvelets \citep{CandesDonohoCurvelets,CandesFDCT05} or bandlets \citep{lepennec-bandelets}, were shown to be very effective in sparsely representing the geometrical content in images. Thanks to the sparsity (or more precisely compressibility) property of these expansions, it is reasonable to assume that essentially only a few large coefficients will contain information about the underlying image, while small values can be attributed to the noise which uniformly contaminates all transform coefficients. Thus, the wavelet thresholding/shrinkage procedure can be mimicked for these transforms, even though some care should be taken when the transform is redundant (corresponding to a frame or a tight frame). The modus operandi is again the same, first apply the transform, then perform a non-linear operator on the coefficients (each coefficient individually or in group of coefficients), and finally apply the inverse transform to get an image estimate. Among the many transform-domain image denoising algorithms to date, we would like to cite \citep{Portilla03,Sendur02,Pizurica02,Luisier07} which are amongst the most efficient in the literature. Except \citep{Luisier07}, all cited approaches use orthodox Bayesian machinery and assume different forms of multivariate priors over blocks of neighboring coefficients and even interscale dependency. Nonetheless, none of those papers provide a study of the theoretical performance of the estimators.

From a theoretical point of view, Cand\`es \citep{Ridgestat} has shown that the ridgelet-based individual coefficient thresholding estimator is nearly minimax for recovering piecewise smooth images away from discontinuities along lines. Individual thresholding of curvelet tight frame coefficients yields an estimator that achieves a nearly-optimal minimax rate $O(n^{-4/3})$\footnote{It is supposed that the image has size $n \times n$.} (up to logarithmic factor) uniformly over the class of piecewise $C^2$ images away from singularities along $C^2$ curves{\textemdash} so-called $C^2$-$C^2$ images \citep{CandesDonohoCurvelets2}\footnote{Known as the cartoon model.}. Similarly, Le Pennec et al. \citep{LePennec07} have recently proved that individual thresholding in an adaptively selected best bandlet orthobasis is nearly-minimax for $C^{\alpha}$ functions away from $C^{\alpha}$ edges.

In the image processing community, block thresholding/shrinkage in a non-Bayesian framework has been used very little. In \citep{Chaux05,Chaux08} the authors propose a multi-channel block denoising algorithm in the wavelet domain. The hyperparameters associated to their method (e.g. threshold), are derived using Stein's risk estimator. Yu et al. \citep{Guo08} advocated the use of BlockJS \citep{cai99} to denoise audio signal in the time-frequency domain with anisotropic block size. To the best of our knowledge, no theoretical study of the minimax properties of block thresholding/shrinkage for images, and more generally for multi-dimensional data, has been reported in the literature.

\subsection{Contributions}
In this paper, we propose a generalization of Stein block thresholding to any dimension $d$. We investigate its minimax properties with a particular emphasis on $d=2$. Towards this goal, we consider a frame coefficient space over which minimaxity is proved; see \eqref{eq:coeffball}. The choice of this space is inspired by the characterization provided in \cite{BorupNielsen} of family of smoothness spaces on $\mathbb{R}^d$, a subclass of so-called decomposition spaces \cite{BorupNielsen,Feichtinger}. We will elaborate more on these (sparsity) smoothness spaces later in subsection~\ref{subsec:stat}. From this characterization, it turns out that our frame coefficient spaces are closely related to smoothness spaces that cover the classical case of Besov spaces, as well as smoothness spaces corresponding to curvelet-type constructions in $\mathbb{R}^d$, $d \geq 2$. Therefore, for $d=2$ our denoiser will apply to both images with smoothness in Besov spaces for which wavelets are known to provide a sparse representation, and also to images that are compressible in the curvelet domain.

Our main theoretical result investigates the minimax rates over these decomposition spaces, and shows that our block estimator can achieve the optimal minimax rate, or is at least nearly-minimax (up to a $\log$ factor) in the least favorable situation. Another novelty is that the minimax rates given here are stated for a general noise sequence model in the transform coefficient domain beyond the usual i.i.d. Gaussian case. Thus, our result is particularly useful when the transform used corresponds to a frame, where a bounded zero-mean white Gaussian noise in the original domain is transformed into a bounded zero-mean correlated Gaussian process with a covariance matrix given by the Gram matrix of the frame.

The choice of the threshold parameter is theoretically discussed and its optimal value is stated for some noise models such as the (non-necessarily i.i.d.) Gaussian case. We provide a simple, fast and a practical procedure. We report a comprehensive simulation study to support our theoretical findings. It turns out that the only two parameters of our Stein block denoiser{\textemdash}the block size and the threshold{\textemdash} dictated by the theory work well for a large set of test images and various transforms. Moreover, the practical performance of our Stein block denoising compares very favorably to state-of-the art methods such as the BLS-GSM of \cite{Portilla03}. Our procedure is however much simpler to implement and has a much lower computational cost than orthodox Bayesian methods such as BLS-GSM, since it does not involve any computationally consuming integration nor optimization steps. A toolbox is made available for download on the Internet to reproduce the results discussed in this paper.

\subsection{Organization of the paper}
The paper is organized as follows. Section \ref{one} is devoted to the one-dimensional BlockJS procedure introduced in \citep{cai99}.
In Section \ref{multi}, we extend BlockJS to the multi-dimensional case and a fairly general noise model beyond the i.i.d. Gaussian case. This section also contains our main theoretical results. In Section \ref{experimental}, a comprehensive experimental study is reported and discussed. We finally conclude in Section \ref{conclusion} and point to some perspectives. The proofs of the results are deferred to the appendix awaiting inspection by the interested reader.

\section{The one-dimensional BlockJS}\label{one}
In this section, we present the construction and the theoretical performance of the one-dimensional BlockJS procedure developed by \citep{cai99}.

Consider the one-dimensional nonparametric regression model:
\begin{equation}\label{mod1D}
Y_i=f(i/n)+ \epsilon_i, \ \ \ \ \ \ i=1, ..., n,
\end{equation}
where $(Y_i)_{i=1, ..., n}$ are the observations, $(\epsilon_i)_{i=1, ..., n}$ are i.i.d. $\mathcal{N}(0,1)$, and $f : [0,1] \rightarrow \mathbb{R}$ is an unknown function. The goal is to estimate $f$ from the observations.
In the orthogonal wavelet framework, \eqref{mod1D} amounts to the sequence model
\begin{equation}\label{model1}
y_{j,k}=\theta_{j, k}+n^{-1/2} z_{j, k}, \ \ \ j=0,...,J, \ \ \ \ \ k=0,...,2^j-1,
\end{equation}
where $J=\lfloor \log_2 n\rfloor$,
$(y_{j,k})_{j,k}$ are the observations, for each $j$, $(z_{j, k})_{k}$ are i.i.d. $\mathcal{N}(0,1)$, and $(\theta_{j,k})_{j,k}$ are approximately the true wavelet coefficients of $f$.
Since they determine completely $f$, the goal is to estimate these coefficients as accurately as possible.
To assess the performance of an estimator $\widehat \theta=(\widehat \theta_{j,k})_{j,k}$ of $\theta=(\theta_{j,k})_{j,k}$, we adopt the minimax approach under the expected squared error over a given Besov body.
The expected squared error is defined by $R(\widehat \theta,\theta) = \sum_{j=0}^{\infty}\sum_{k=0}^{2^j-1}\mathbb{E}\left((\widehat \theta_{j,k}-\theta_{j,k})^2\right)$,
and the Besov body by
$${\it \Theta}_{p,q}^s(M) =\left\lbrace \theta=(\theta_{j, k})_{j, k} ; \ \ \ \ \left( \sum_{j=0}^{\infty}\left(2^{j(s+1/2-1/p)}\left( \sum_{k =0}^{2^j-1}|\theta_{j,k}|^p\right)^{1/p} \right)^q \right)^{1/q} \le M  \right\rbrace.$$
In this notation, $s>0$ is a smoothness parameter, $0 < p \le +\infty$ and $0 < q \le +\infty$ are norm parameters\footnote{This is a slight abuse of terminology as for $0 < p,q < 1$, Besov spaces are rather complete quasinormed linear spaces.}, and $M\in (0,\infty)$ denotes the radius of the ball. The Besov body contains a wide class of $\theta=(\theta_{j, k})_{j, k}$. It includes the H\"older body ${\it \Theta}^s_{\infty,\infty}(M)$ and the Sobolev body ${\it \Theta}^s_{2,2}(M)$.

The goal of the minimax approach is to construct an adaptive estimator $\widehat \theta=(\widehat \theta_{j,k})_{j,k}$ such that $\sup_{\theta \in {\it \Theta}_{p,q}^s(M)}R(\widehat \theta,\theta)$ is as small as possible.
A candidate is the BlockJS procedure whose paradigm is described below.

Let $L= \lfloor \log n \rfloor$ be the block size, $j_0=\lfloor  \log_2 L \rfloor$ the coarsest decomposition scale and, for any $j$, $\Lambda_j =\{0,...,2^j-1\}$ is the set of locations at scale $j$. For any $j\in\{j_0,...,J\}$, let ${A}_j=\lbrace 1,..., \lfloor 2^{j}L^{-1}\rfloor \rbrace$ be the set of block indices at scale $j$, and for any $K \in {A}_j$, $U_{j,K} = \lbrace k \in \Lambda_j; \, (K-1)L\le k\le K L-1 \rbrace$ is the set indexing the locations of coefficients within the $K$th block.
Let $\lambda_*$ be a threshold parameter chosen as the root of $x-\log x=3$ (i.e. $\lambda_*= 4.50524...$). Now estimate $\theta=(\theta_{j,k})_{j,k}$ by $\widehat \theta^*=(\widehat \theta^*_{j,k})_{j,k}$ where, for any $k \in U_{j,K}$ and $K \in {A}_j$,
\begin{eqnarray}\label{seii}
\widehat \theta^*_{j,k}=
\begin{cases}
y_{j,k},  & {\text{if $j \in \{0,...,j_0-1\}$}} ,\\
y_{j,k} \left(1- \frac{\lambda_* n^{-1}}{\tfrac{1}{L}\sum_{k\in
U_{j,K}}y_{j,k}^2}\right)_{+}, & {\text{if  $j \in \{j_0,...,J\}$}}, \\
0, & {\text{if $j\in \mathbb{N}-\{0,...,J\}$}}.
\end{cases}
\end{eqnarray}
where $(x)_+ = \max(x,0)$. Thus, at the coarsest scales $j \in \{0,...,j_0\}$, the observed coefficients $\left(y_{j,k}\right)_{k \in \Lambda_j}$ are left intact as usual. For $k\in \Lambda_j$ and $j\in \mathbb{N}-\{0,...,J\}$, $\theta_{j,k}$ is estimated by zero. For $k \in U_{j,K}$, $K \in {A}_j$ and $j \in \{j_0,...,J\}$, if the mean energy within the $K$th block $\sum_{k\in U_{j,K}}y_{j,k}^2/L$ is larger than $\lambda_* n^{-1}$ then $y_{j,k}$ is shrunk by the amount $y_{j,k}\frac{\lambda_* n^{-1}}{\frac{1}{L}\sum_{k\in U_{j,K}}y_{j,k}^2}$; otherwise, $\theta_{j,k}$ is estimated by zero. Note that $\frac{\frac{1}{L}\sum_{k\in U_{j,K}}y_{j,k}^2}{n^{-1}}$ can be interpreted as a local measure of signal-to-noise ratio in the block $U_{j,K}$. Such a block thresholding originates from the James-Stein rule introduced in \citep{stein}.

The block length $L= \lfloor \log n \rfloor$ and the value $\lambda_*=4.50524$ are chosen based on theoretical considerations;
under this calibration, the BlockJS is (near) optimal in terms of minimax rate and adaptivity. This is summarized in the following theorem.
\begin{theorem}[\citep{cai99}]\label{mini}
Consider the model \eqref{model1} for $n$ large enough. Let $\widehat \theta^*$ be given as \eqref{seii}. Then there exists a constant $C>0$ such that
\begin{eqnarray}\label{rea}
\sup_{\theta \in {\it \Theta}_{p,q}^s (M)} R(\widehat \theta^*,\theta)\le C
\left\{
\begin{aligned}
& n^{-2s/(2s+1)}, & & {\text{for $p\ge 2$}},\\
& n^{-2s/(2s+1)} (\log n)^{(2-p)/(p(2s+1))} ,&  & {\text{for $p<2$, $sp\ge 1$}}.
\end{aligned}
\right.
\end{eqnarray}
\end{theorem}
The rates of convergence \eqref{rea} are optimal, except in the case $p<2$ where there is an extra logarithmic term.
They are better than those achieved by standard individual thresholding (hard, soft, non-negative garotte, etc); we gain a logarithmic factor for $p\ge 2$. See \citep{donohoj}.

\section{The multi-dimensional BlockJS}\label{multi}
This section is the core of our proposal where we introduce a BlockJS-type procedure for multi-dimensional data.
The goal is to adapt its construction in such a way that it preserves its optimal properties over a wide class of functions.

\subsection{The sequence model}
\label{subsec:seq}
Our approach begins by projecting the model \eqref{mod} onto a collection of atoms $\left(\varphi_{j,\ell,{\bf k}}\right)_{j,\ell,{\bf k}}$ that forms a (tight) frame. This gives rise to a sequence space model obtained by calculating the noisy coefficient $y_{j,\ell,{\bf k}}=\langle Y,\varphi_{j,\ell,{\bf k}}\rangle$ for any element of the frame $\varphi_{j,\ell,{\bf k}}$. We then have a multi-dimensional sequence of coefficients $(y_{j,\ell,{\bf k}})_{j,\ell,{\bf k}}$ defined by
\begin{equation}\label{model}
y_{j,\ell,{\bf k}}=\theta_{j,\ell,{\bf k}}+n^{-r/2}z_{j,\ell,{\bf k}}, \ \ \ j =0,...,J, \ \ \ \ \ell \in B_j, \ \ \ \ {\bf k}\in D_j,
\end{equation}
where $J=\lfloor \log_2 n \rfloor$, $r\in [1,d]$, $d\in \mathbb{N}^*$, $B_j=\{1,..., \lfloor c_* 2^{\upsilon j}\rfloor\}$, $c_*\ge 1$, $\upsilon \in [0,1]$, ${\bf k}=(k_1,...,k_d)$, $D_j=\prod_{i=1}^d\{0,..., \lfloor 2^{\mu_i j} \rfloor -1\}$, $(\mu_i)_{i =1,...,d}$ is a sequence of positive real numbers, $(z_{j,\ell,{\bf k}})_{j,\ell,{\bf k}}$ are random variables and $(\theta_{j,\ell,{\bf k}})_{j,\ell,{\bf k}}$ are unknown coefficients. Let $d_*=\sum_{i=1}^d \mu_i$.

The indices $j$ and $\bf{k}$ are respectively the scale and position parameters. $\ell$ is a generic integer indexing for example the orientation (subband) which may be scale-dependent. The parameters $(\mu_i)_{i =1,...,d}$ allow to handle anisotropic subbands. To illustrate the meaning of these parameters, let's see how they specialize in some popular transforms. For example, with the separable two-dimensional wavelet transform, we have $v=0$, $c^*=3$, and $\mu_1=\mu_2=1$. Thus, as expected, we get three isotropic subbands at each scale. For the second generation curvelet transform \cite{CandesFDCT05}, we have $v=1/2$, $\mu_1=1$ and $\mu_2=1/2$ which corresponds to the parabolic scaling of curvelets.

\subsubsection{Assumptions on the noise sequence}
\label{subsubsec:noiseseq}
Let $L= \lfloor ( r \log n)^{1/d} \rfloor$ be the block length, $j_0=\left\lfloor  (1/\min_{i=1,...,d} \mu_i) \log_2 L \right\rfloor$ is the coarsest decomposition scale, and $J_*=\lfloor (r/(d_*+\delta+\upsilon))\log_2 n\rfloor$.
For any $j\in\{j_0,...,J_*\}$, let
\begin{itemize}
\item $\mathcal{A}_j=\prod_{i=1}^d \lbrace 1,..., \lfloor 2^{\mu_i j}L^{-1}\rfloor \rbrace$ be the set indexing the blocks at scale $j$.
\item For each block index ${\bf K}=(K_1,...,K_d) \in \mathcal{A}_j$, $U_{j,{\bf K}} = \lbrace {\bf k} \in D_j; \, (K_1-1)L\le k_1\le K_1 L-1,..., \, (K_d-1)L\le k_d\le K_d L-1 \rbrace$ is the set indexing the positions of coefficients within the ${\bf K}$th block $U_{j,{\bf K}}$.
\end{itemize}
Our assumptions on the noise model are as follows. Suppose that there exist $\delta \ge 0$, $\lambda_*>0$, $Q_1>0$ and $Q_2>0$ independent of $n$ such that
\begin{enumerate}[$(\mathrm{A}1)$]
\item $\sup_{j\in \{0,...,J\}}\sup_{\ell\in B_j}2^{-j(d_*+\delta)} \sum_{{\bf k} \in D_j} \mathbb{E}\left( z^2_{j,\ell,{\bf k}}\right)\le Q_1$.
\item
\begin{equation*}
\sum_{j=j_0}^{J_*}\sum_{\ell \in B_j}\sum_{{\bf K} \in \mathcal{A}_j}  \sum_{{\bf k} \in U_{j,{\bf K}} } \mathbb{E}\left( z_{j,\ell,{\bf k}}^2\indic{\left\lbrace  \sum_{{\bf k} \in U_{j,{\bf K}} } z_{j,\ell,{\bf k}}^2 > \lambda_*2^{\delta j}L^d/4 \right\rbrace}\right)  \le Q_2.
\end{equation*}
\end{enumerate}

Assumptions $(\mathrm{A}1)$ and $(\mathrm{A}2)$ are satisfied for a wide class of noise models on the sequence $(z_{j,\ell,{\bf k}})_{j,\ell,{\bf k}}$ (not necessarily independent or identically distributed). Several such noise models are characterized in Propositions \ref{mich} and \ref{mich2} below.

\begin{remark}[\it Comments on $\delta$]
The parameter $\delta$ is connected to the nature of the model.
For standard models, and in particular, the $d$-dimensional nonparametric regression corresponding to the problem of denoising (see Section \ref{experimental}), $\delta$ is set to zero.
The presence of $\delta$ in our assumptions, definitions and results is motivated by potential applicability of the multi-dimensional BlockJS (to be defined in Subsection \ref{dd}) to other inverse problems such as deconvolution. 
The role of $\delta$ becomes of interest when addressing such inverse problems. This will be the focus of a future work. To illustrate the importance of $\delta$ in one-dimensional deconvolution, see \citep{johnstonek}.
\end{remark}

\subsection{The smoothness space}
\label{subsec:stat}

We wish to estimate $(\theta_{j,\ell,{\bf k}})_{j,\ell,{\bf k}}$ from $(y_{j,\ell,{\bf k}})_{j,\ell,{\bf k}}$ defined by \eqref{model}. To measure the performance of an estimator $\widehat \theta=(\widehat \theta_{j,\ell,{\bf k}})_{j,\ell,{\bf k}}$ of $\theta=(\theta_{j,\ell,{\bf k}})_{j,\ell,{\bf k}}$, we consider the minimax approach under the expected multi-dimensional squared error over a multi-dimensional frame coefficient space.
The expected multi-dimensional squared error is defined by $$R\left(\widehat \theta,\theta\right) = \sum_{j=0}^{\infty}\sum_{\ell\in B_j}\sum_{{\bf k} \in D_j}\mathbb{E}\left((\widehat \theta_{j,\ell,{\bf k}}-\theta_{j,\ell,{\bf k}})^2\right)$$ and the multi-dimensional frame coefficient smoothness/sparseness space by
\begin{equation}
\label{eq:coeffball}
{\boldsymbol\Theta}_{p,q}^s(M) =\left\lbrace \theta=(\theta_{j,\ell,{\bf k}})_{j,\ell,{\bf k}} ; \ \ \left( \sum_{j=0}^{\infty}\sum_{\ell\in B_j} \left(2^{j(s+d_*/2-d_*/p)}\left( \sum_{{\bf k} \in D_j}|\theta_{j,\ell,{\bf k}}|^p\right)^{1/p} \right)^q \right)^{1/q} \le M  \right\rbrace,
\end{equation}
with a smoothness parameter $s$, $0 < p \le +\infty$ and $0 < q \le +\infty$. We recall that $d_*=\sum_{i=1}^d \mu_i$.

The definition of these smoothness spaces is motivated by the work of \cite{BorupNielsen}. These authors studied decomposition spaces associated to appropriate structured uniform partition of the unity in the frequency space $\mathbb{R}^d$. They considered construction of tight frames adapted to form atomic decomposition of the associated decomposition spaces, and established norm equivalence between these smoothness/sparseness spaces and the sequence norm defined in \eqref{eq:coeffball}. That is, the decomposition space norm can be completely characterized by the sparsity or decay behavior of the associated frame coefficients.

For example, in the case of a "uniform" dyadic partition of the unity, the smoothness/sparseness space is a Besov space $B^{s}_{p,q}$, for which suitable wavelet expansion\footnote{With a wavelet having sufficient regularity and number of vanishing moments \cite{Meyer92}.} is known to provide a sparse representation \cite{Meyer92}. In this case, from subsection \ref{subsec:seq} we have $d^*=d$, and ${\boldsymbol\Theta}_{p,q}^s(M)$ is a $d$-dimensional Besov ball.

Curvelets in arbitrary dimensions correspond to partitioning the frequency plane into dyadic coronae, which are then angularly localized near regions of side length $2^j$ in the radial direction and $2^{j/2}$ in all the other directions \cite{CandesDemanet}. For $d=2$, the angular wedges obey the parabolic scaling law \cite{CandesDonohoCurvelets}. This partition of the frequency plane is significantly different from dyadic decompositions, and as a consequence, sparseness for curvelet expansions cannot be described in terms of classical smoothness spaces. For $d=2$, Borup and Nielsen \cite[Lemma 10]{BorupNielsen} showed that the smoothness/sparseness space \eqref{eq:coeffball} and the smoothness/sparseness of the second-generation curvelets \cite{CandesFDCT05} are the same, in which case $d^*=3/2$. Embedding results for curvelet-type decomposition spaces relative to Besov spaces were also provided in \cite{BorupNielsen}. Furthermore, it was shown that piecewise $C^2$ images away from piecewise-$C^2$ singularities, which are sparsely represented in the curvelet tight frame \cite{CandesDonohoCurvelets}, are contained in ${\boldsymbol\Theta}_{2/3,2/3}^{3/2 + \beta}$, $\forall \beta>0$. Even though the role and the range of $\beta$ has not been clarified by the authors in \cite{BorupNielsen}.

\subsection{Multi-dimensional block estimator}\label{dd}
As for the one-dimensional case, we wish to construct an adaptive estimator $\widehat \theta=(\widehat \theta_{j,\ell,{\bf k}})_{j,\ell,{\bf k}}$ such that $\sup_{\theta \in {\boldsymbol\Theta}_{p,q}^s(M)}R\left(\widehat \theta,\theta\right)$ is as small as possible. To reach this goal, we propose a multi-dimensional version of the BlockJS procedure introduced in \citep{cai99}.

From subsection \ref{subsubsec:noiseseq}, recall the definitions of $L$, $j_0$, $J_*$, $\mathcal{A}_j$ and $U_{j,{\bf K}}$.
We estimate $\theta=(\theta_{j,\ell,{\bf k}})_{j,\ell,{\bf k}}$ by $\widehat \theta^*=(\widehat \theta^*_{j,\ell,{\bf k}})_{j,\ell,{\bf k}}$ where, for any ${\bf k} \in U_{j,{\bf K}}$, ${\bf K} \in \mathcal{A}_j$ and $\ell\in B_j$,
\begin{eqnarray}\label{sei}
\widehat \theta^*_{j,\ell,{\bf k}}=
\begin{cases}
y_{j,\ell,{\bf k}},  & {\text{if  $j \in \{0,...,j_0-1\}$}} ,\\
y_{j,\ell,{\bf k}} \left(1- \frac{\lambda_* n^{-r} 2^{\delta j}}{\tfrac{1}{L^d}\sum_{{\bf k}\in
U_{j,{\bf K}}}y_{j,\ell,{\bf k}}^2}\right)_{+}, & {\text{if $j \in \{j_0,...,J_*\}$}}, \\
0 , & {\text{if $j\in \mathbb{N}-\{0,...,J_*\}$}}.
\end{cases}
\end{eqnarray}
In this definition, $\delta$ and $\lambda_*$ denote the constants involved in $(\mathrm{A}1)$ and $(\mathrm{A}2)$.
Again, the coarsest scale coefficients are left unaltered, while the other coefficients are either thresholded or shrunk depending whether the local measure of signal-to-noise ratio $\frac{\frac{1}{L^d}\sum_{{\bf k}\in
U_{j,{\bf K}}}y_{j,\ell,{\bf k}}^2}{n^{-r}}$ within the block $U_{j,{\bf K}}$ is larger that the threshold $\lambda_* 2^{\delta j}$.
Notice that the dimension $d$ of the model appears in the definition of $L$, the length of each block $U_{j,{\bf K}}$.
This point is crucial; $L$ optimizes the theoretical and practical performance of the considered multi-dimensional BlockJS procedure.
As far as the choice of the threshold parameter $\lambda_*$ is concerned, it will be discussed in Subsection \ref{lambda} below.

\subsection{Minimax theorem}
Theorem \ref{mini2} below investigates the minimax rate of \eqref{sei} over ${\boldsymbol\Theta}_{p,q}^s$.
\begin{theorem}\label{mini2}
Consider the model \eqref{model} for $n$ large enough. Suppose that $(\mathrm{A}1)$ and $(\mathrm{A}2)$ are satisfied. Let $\widehat \theta^*$ be given as in \eqref{sei}.
\begin{itemize}
\item There exists a constant $C>0$ such that
\begin{equation*}
\sup_{\theta \in {\boldsymbol\Theta}_{p,q}^s (M)}R\left(\widehat \theta^*,\theta\right) \le C \rho_n,
\end{equation*}
where
\begin{eqnarray}\label{rates}
\rho_n =  \left\{
\begin{aligned}
& n^{-2sr/(2s+\delta+d_*+\upsilon)}, & & {\text{for $q\le 2 \le p$}},\\
& (\log n/n)^{2sr/(2s+\delta+d_*+\upsilon)}, &  & {\text{for $q \le p<2$, $sp>d_*\vee (1-p/2)(\delta+d_*+\upsilon)$}}.
\end{aligned}
\right.
\end{eqnarray}
\item If $\upsilon =0$, the minimax rates \eqref{rates} hold without the restriction $q\le p\wedge 2$.
\end{itemize}
\end{theorem}

The rates of convergence \eqref{rates} are optimal for a wide class of variables $(z_{j,\ell,{\bf k}})_{j,\ell,{\bf k}}$.
If we take $d_*=d=\mu_1=1$, $r=1$, $c_*=1$ and $\upsilon=\delta=0$, then
we recover the rates exhibited in the one-dimensional wavelet case expressed in Theorem \ref{mini}.
There is only a minor difference on the power of the logarithmic term for $p<2$.
Thus, Theorem \ref{mini2} can be viewed as a generalization of Theorem \ref{mini}.

In the case of $d$-dimensional isotropic Besov spaces, where wavelets (corresponding to $\upsilon=0$, $\mu_1=\mu_2=1$ and then $d_*=d$) provide optimally sparse representations, Theorem \ref{mini2} can be applied without the restriction $q \le p\wedge 2$. Therefore, for $p \geq 2$, Theorem \ref{mini2} states that Stein block thresholding gets rid of the logarithmic factor, hence achieving the optimal minimax rate over those Besov spaces. For $p < 2$, the block estimator is nearly-minimax.

As far as curvelet-type decomposition spaces are concerned, from section \ref{subsec:seq} we have $\mu_1=1, \mu_2=\frac{1}{2}, d_*=\mu_1+\mu_2=\frac{3}{2}$, $r=d=2$, $\upsilon =\frac{1}{2}$, $\delta=0$. This gives the rates
\begin{eqnarray*}
\rho_n =  \left\{
\begin{aligned}
& n^{-2s/(s+1)}, & & {\text{for $q \le 2 \le p$}},\\
& (\log n/n)^{2s/(s+1)}, &  & {\text{for $q \le p < 2$, $sp>\frac{3}{2}\vee (2-p)$}}.
\end{aligned}
\right.
\end{eqnarray*}
where the logarithmic factor disappears only for $q \leq 2 \leq p$. Following the discussion of section \ref{subsec:stat}, $C^2$-$C^2$ images correspond to a smoothness space ${\boldsymbol\Theta}_{p,q}^{s}$ with $p=q=2/3$. Moreover, $\exists \kappa > 0$ such that taking $s=2+\kappa$ satisfies the condition of Theorem \ref{mini2}, and $C^2$-$C^2$ images are contained in ${\boldsymbol\Theta}_{2/3,2/3}^{s}$ with such a choice. We then arrive at the rate $O(n^{-4/3})$ (ignoring the logarithmic factor). This is consistent with the results of \cite{KorostelevTsybakov}, which established that no estimator can achieve a better rate than the optimal minimax rate $O(n^{-4/3})$ uniformly over the $C^2$-$C^2$ class. On the other hand, individual thresholding in the curvelet tight frame has also the nearly-minimax rate $O(n^{-4/3})$ \cite{CandesDonohoCurvelets2} uniformly over the class of $C^2$-$C^2$ images. Nonetheless, the experimental results reported in this paper indicate that block curvelet thresholding outperforms in practice term-by-term thresholding on a wide variety of images, although the improvement can be of a limited extent.



\subsection{On the (theoretical) choice of the threshold}\label{lambda}

To apply Theorem \ref{mini2}, it is enough to determine $\delta$ and $\lambda_*$ such that $(\mathrm{A}1)$ and $(\mathrm{A}2)$ are satisfied.
The parameter $\delta$ is imposed by the nature of the model; it can be easily fixed as in our denoising experiments where it was set to $\delta=0$. The choice of the threshold $\lambda_*$ is more involved. This choice is crucial towards good performance of the estimator $\widehat \theta^*$.
From a theoretical point of view, since the constant $C$ of the bound \eqref{rates} increases with growing $\lambda$, the optimal threshold is the smallest real number $\lambda_*$ such that $(\mathrm{A}2)$ is fulfilled. In the following, we first provide the explicit expression of $\lambda_*$ in the situation of a non-necessarily i.i.d. Gaussian noise sequence $(z_{j,\ell,{\bf k}})_{j,\ell,{\bf k}}$. This result is then refined in the case of a white Gaussian noise.

Proposition \ref{mich} below determines a suitable threshold $\lambda_*$ satisfying $(\mathrm{A}1)$ and $(\mathrm{A}2)$ when $(z_{j,\ell,{\bf k}})_{j,\ell,{\bf k}}$ are Gaussian random variables (not necessarily i.i.d.).
\begin{proposition}\label{mich}
Consider the model \eqref{model} for $n$ large enough. Suppose that, for any $j \in \{0,...,J\}$ and any $\ell\in B_j$, $(z_{j,\ell,{\bf k}})_{{\bf k}}$ is a centered Gaussian process. Assume that there exists two constants $Q_3>0$ and $Q_4>0$ (independent of $n$) such that
\begin{itemize}
\item $(\mathrm{A}3)$: $\sup_{j\in \{0,...,J\}}\sup_{\ell \in B_j} \sup_{{\bf k}\in D_j} 2^{-2\delta j} \mathbb{E}\left(z_{j,\ell,{\bf k}}^4\right)\le Q_3$.
\item $(\mathrm{A}4)$: for any $a=(a_{\bf k})_{{\bf k} \in D_j}$ such that $\sup_{j\in \{0,...,J\}} \sup_{{\bf K} \in \mathcal{A}_j} \sum_{{\bf k}\in U_{j,{\bf K}}}a_{\bf k}^{2}\le  1$, we have
$$\sup_{j\in \{0,...,J\}}\sup_{\ell\in B_j}\sup_{{\bf K} \in \mathcal{A}_j} 2^{-\delta j}\mathbb{E}\left( \left( \sum_{{\bf k}\in U_{j,{\bf K}}}  a_{{\bf k}}z_{j,\ell,{\bf k}}\right)^2\right)  \le Q_4.$$
\end{itemize}
Then $(\mathrm{A}1)$ and $(\mathrm{A}2)$ are satisfied with $\lambda_*=4\left( (2Q_4)^{1/2}+Q_3^{1/4}\right)^2$. Therefore Theorem \ref{mini2} can be applied to
$\widehat \theta^*$ defined by \eqref{sei} with such a $\lambda_*$.
\end{proposition}
This result is useful as it establishes that the block denoising procedure and the minimax rates of Theorem \ref{mini} apply to the case of frames where a bounded zero-mean white Gaussian noise in the original domain is transformed into a bounded zero-mean correlated Gaussian process.

If additional information is considered on $(z_{j,\ell,{\bf k}})_{j,\ell,{\bf k}}$, the threshold constant $\lambda_*$ defined in Proposition \ref{mich} can be improved. This is the case when $(z_{j,\ell,{\bf k}})_{j,\ell,{\bf k}}$ are i.i.d. $\mathcal{N}(0,1)$ as is the case if the transform were orthogonal (e.g. orthogonal wavelet transform). The statement is made formal in the following proposition.
\begin{proposition}\label{mich2}
Consider the model \eqref{model} for $n$ large enough. Suppose that, for any $j \in \{0,...,J\}$ and any $\ell\in B_j$, $(z_{j,\ell,{\bf k}})_{{\bf k}}$ are i.i.d. $\mathcal{N}(0,1)$ as is the case when the transform used corresponds to an orthobasis.
Theorem \ref{mini2} can be applied with the estimator
$\widehat \theta^*$ defined by \eqref{sei} with $\delta=0$ and $\lambda_*$ the root of $x-\log x=3$, i.e. $\lambda_*= 4.50524...$ .
\end{proposition}
The optimal threshold constant $\lambda_*$ described in Proposition \ref{mich2} corresponds to the one isolated by \citep{cai99}.

\section{Application to image block denoising}\label{experimental}
\subsection{Impact of threshold and block size}
\begin{figure}
\centerline{Barbara $512 \times 512$}
\psfig{figure=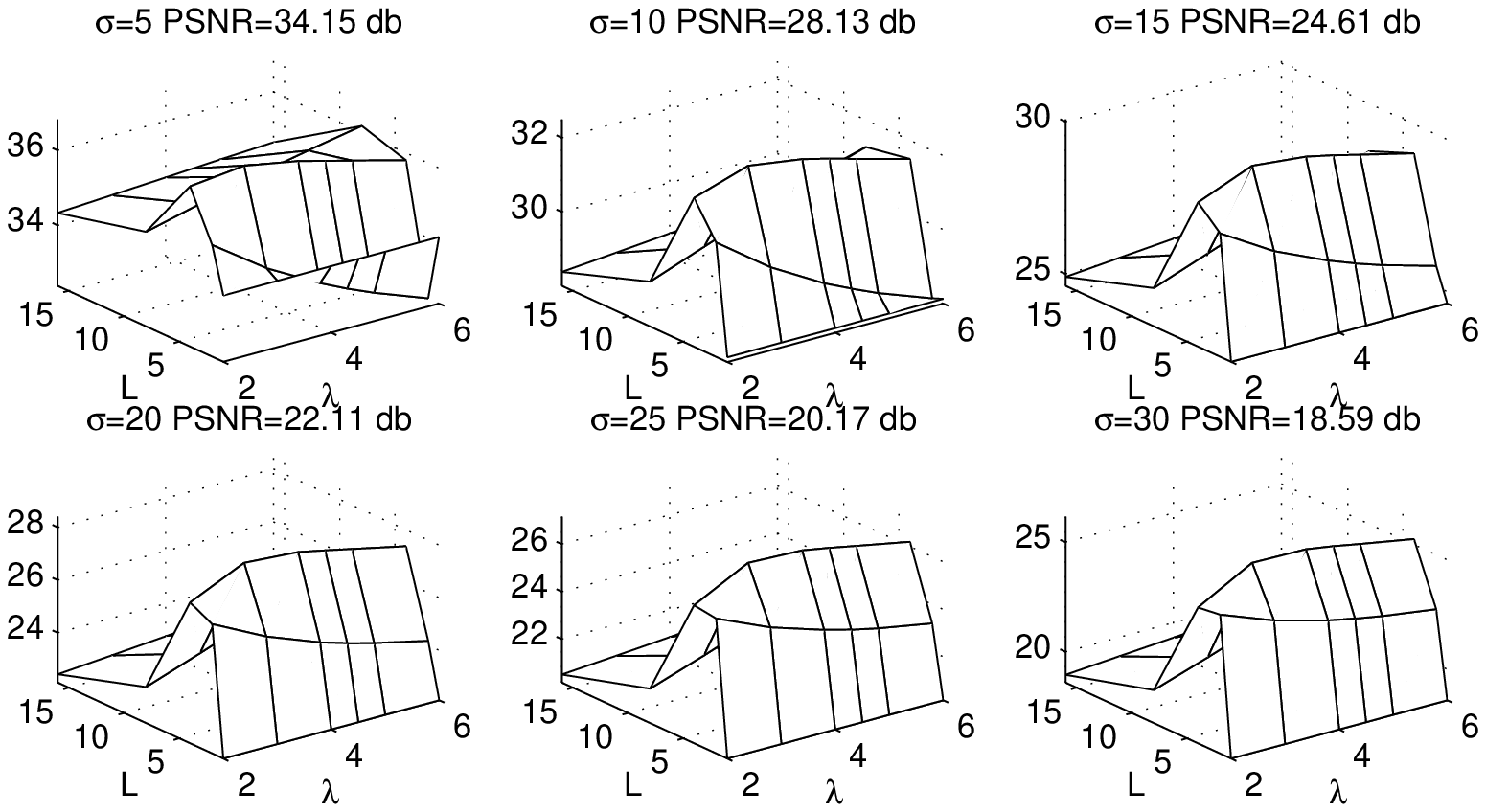,bbllx=3cm,bblly=12cm,bburx=19cm,bbury=21cm,width=\textwidth,clip=}
\centerline{Peppers $256 \times 256$}
\psfig{figure=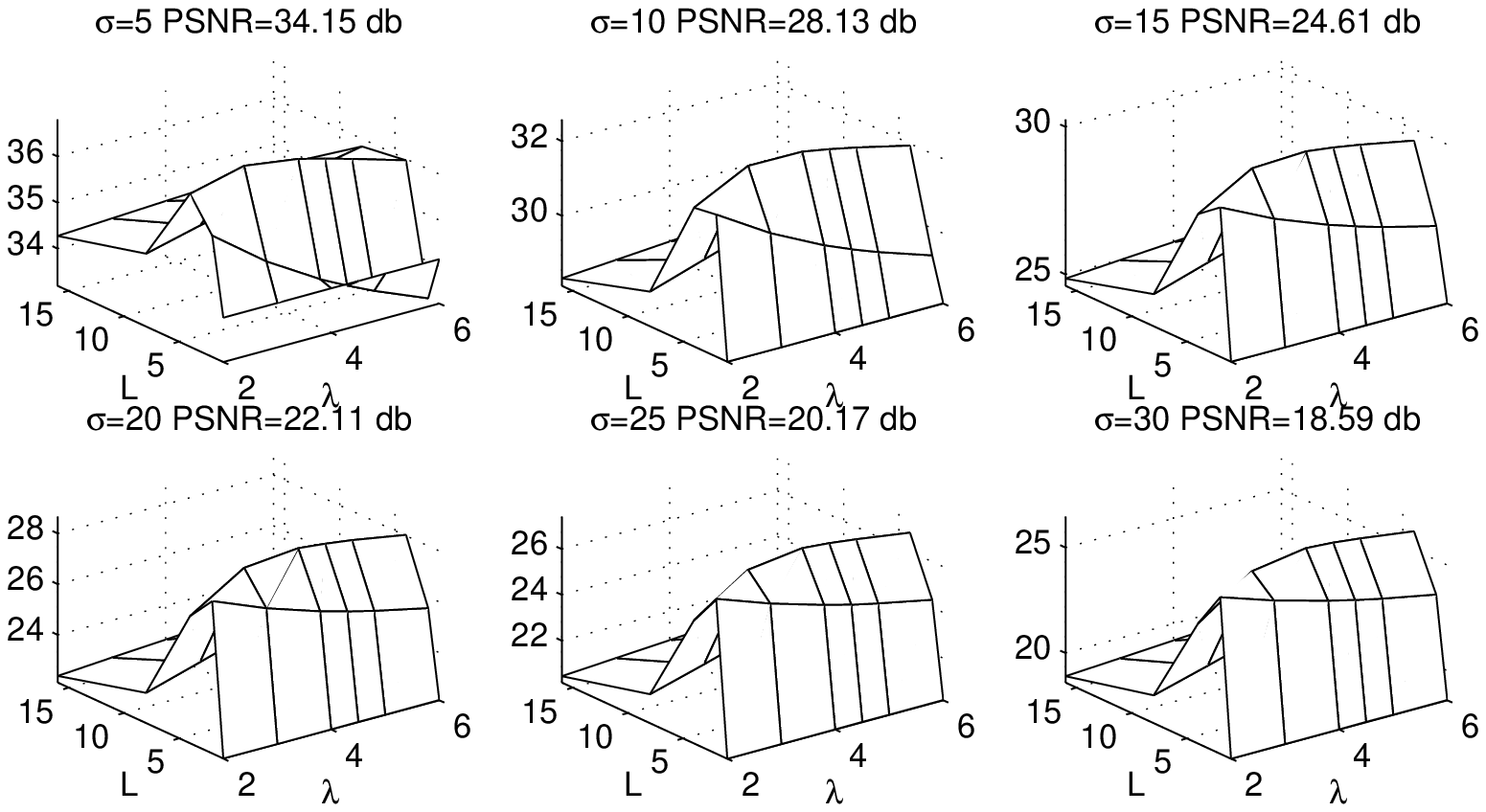,bbllx=3cm,bblly=12cm,bburx=19cm,bbury=21cm,width=\textwidth,clip=}
\caption{Output PSNR as a function of the block size and the threshold $\lambda$ at different noise levels $\sigma \in \{5,10,15,20,25,30\}$. Block denoising was applied in the DWT domain.}
\label{fig:dwtimpact}
\end{figure}

\begin{figure}
\centerline{Barbara $512 \times 512$}
\psfig{figure=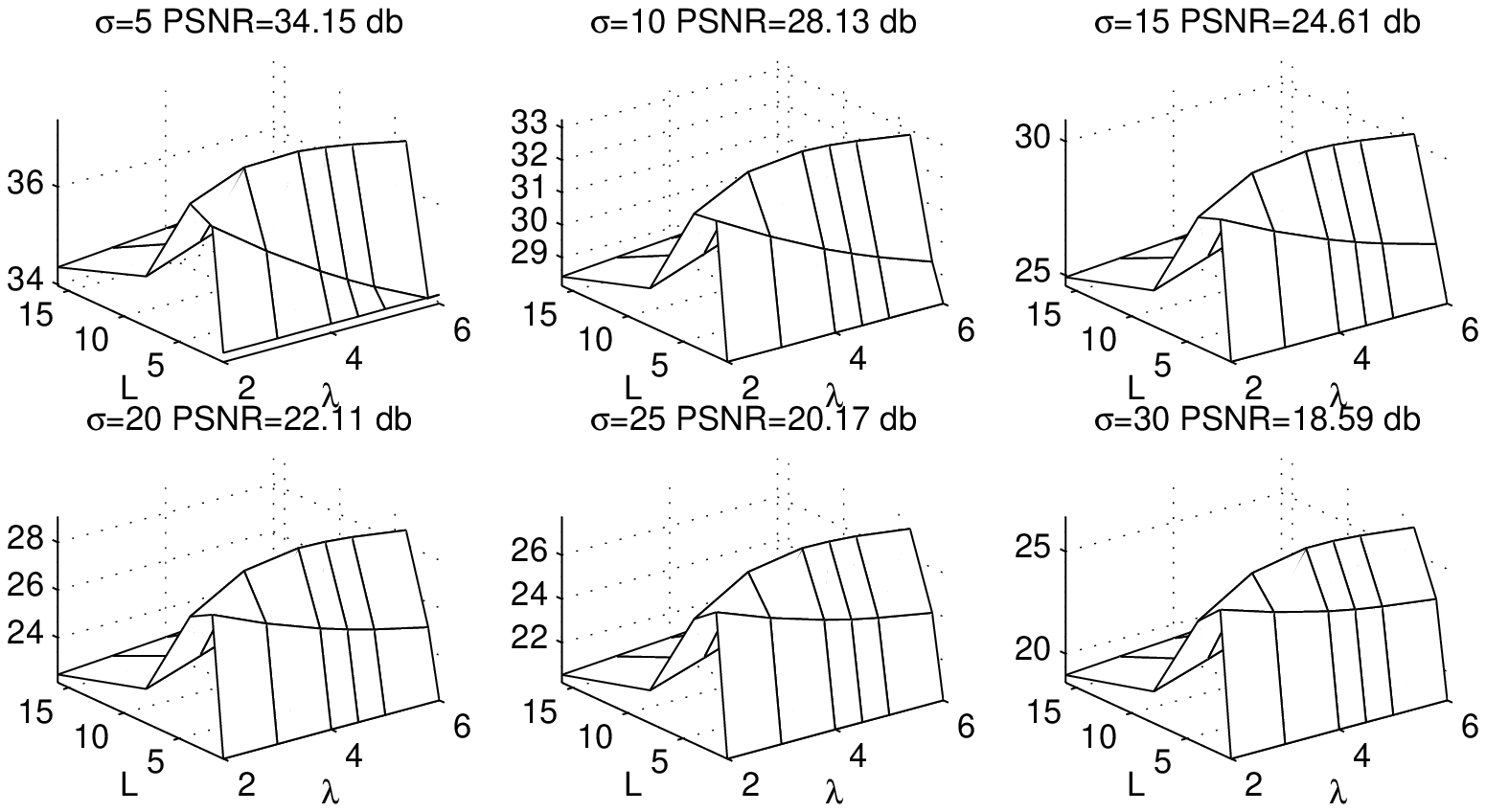,bbllx=3cm,bblly=12cm,bburx=19cm,bbury=21cm,width=\textwidth,clip=}
\centerline{Peppers $256 \times 256$}
\psfig{figure=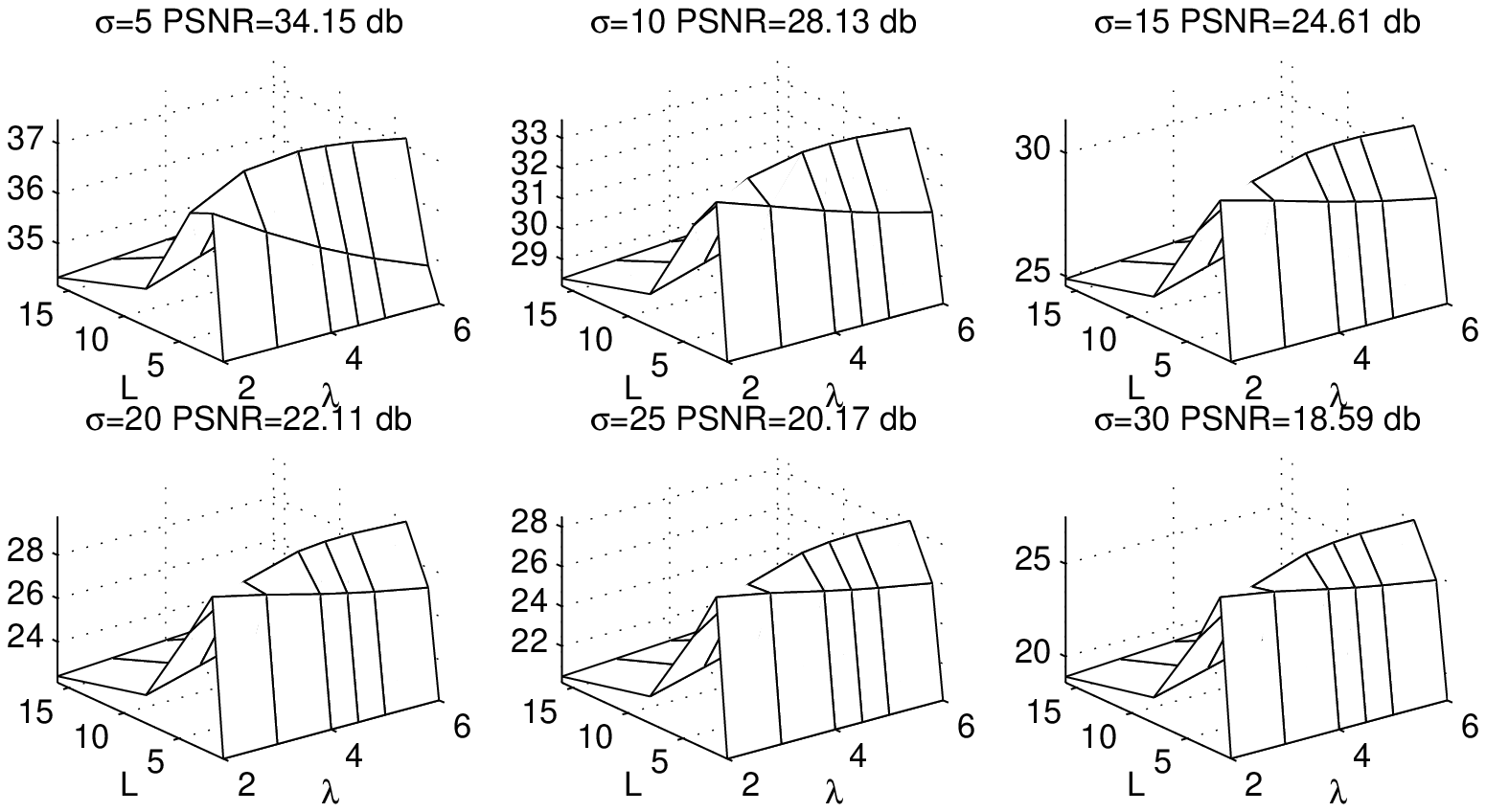,bbllx=3cm,bblly=12cm,bburx=19cm,bbury=21cm,width=\textwidth,clip=}
\caption{Output PSNR as a function of the block size and the threshold $\lambda$ at different noise levels $\sigma \in \{5,10,15,20,25,30\}$. Block denoising was applied in the UDWT domain.}
\label{fig:udwtimpact}
\end{figure}

\begin{figure}
\centerline{Barbara $512 \times 512$}
\psfig{figure=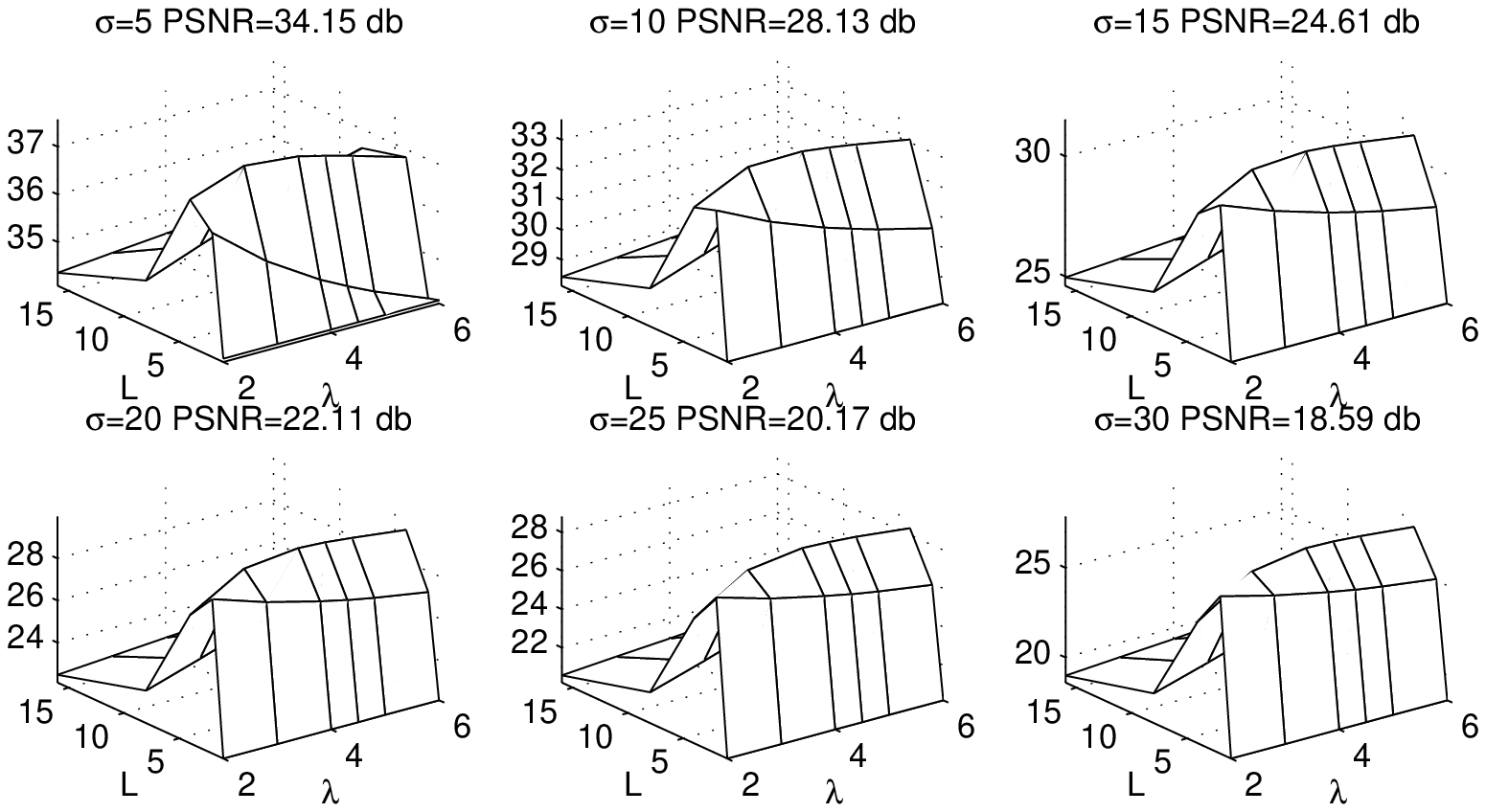,bbllx=3cm,bblly=12cm,bburx=19cm,bbury=21cm,width=\textwidth,clip=}
\centerline{Peppers $256 \times 256$}
\psfig{figure=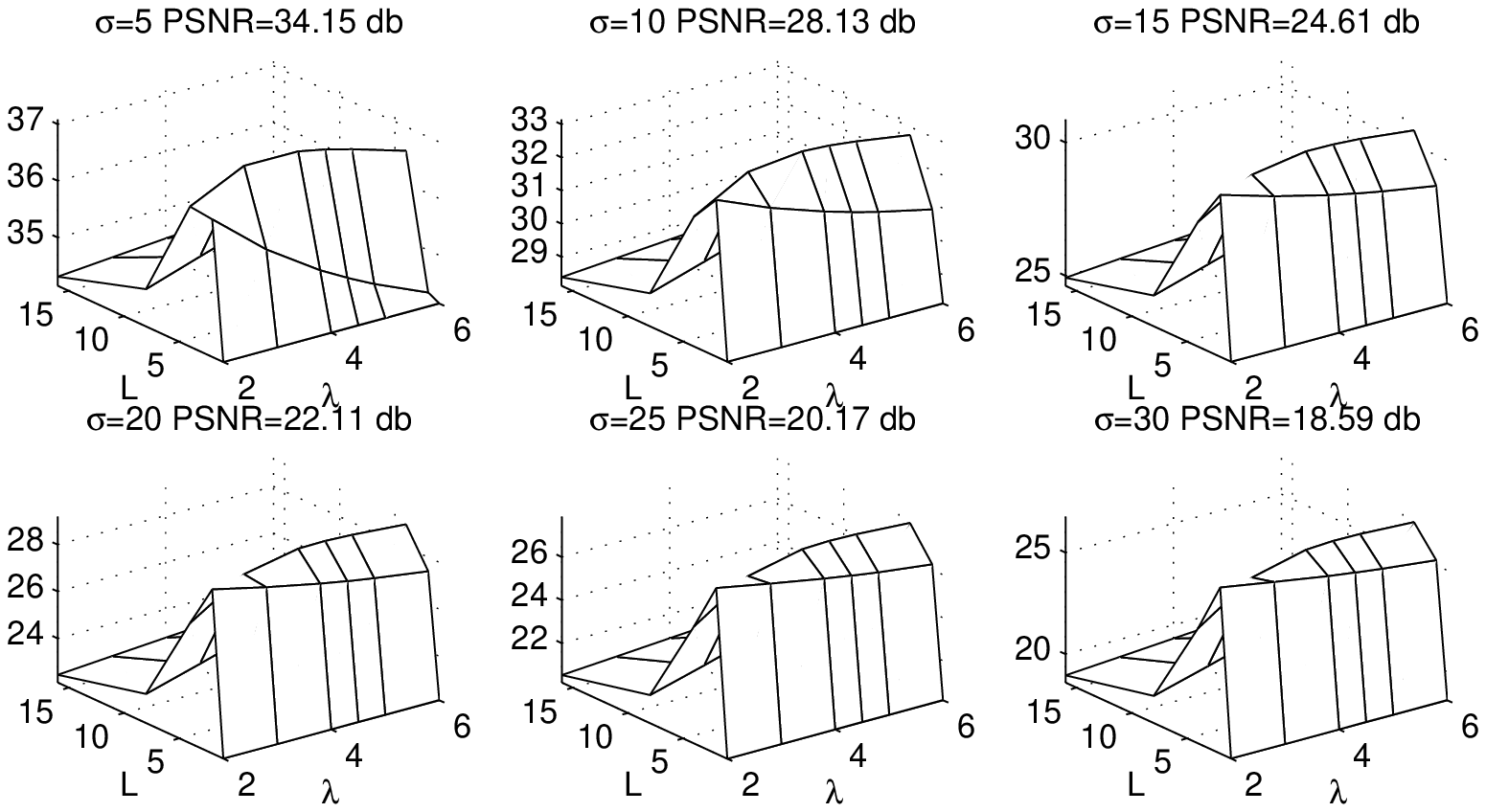,bbllx=3cm,bblly=12cm,bburx=19cm,bbury=21cm,width=\textwidth,clip=}
\caption{Output PSNR as a function of the block size and the threshold $\lambda$ at different noise levels $\sigma \in \{5,10,15,20,25,30\}$. Block denoising was applied in the FDCT domain.}
\label{fig:curvimpact}
\end{figure}

In this first experiment, the goal is twofold: first assess the impact of the threshold and the block size on the performance of block denoising, and second investigate the validity of their choice as prescribed by the theory. For a $n \times n$ image $f$ and its estimate $\hat{f}$, the denoising performance is measured in terms of peak signal-to-noise ratio (PSNR) in decibels (dB)
\[
\psnr = 20 \log_{10} \frac{n\|f\|_\infty}{\|\hat{f}-f\|_2} ~ \mathrm{dB} ~.
\]
In this experiment, as well as in the rest of paper, three popular transforms are used: the orthogonal wavelet transform (DWT), its translation invariant version (UDWT) and the second generation fast discrete curvelet transform (FDCT) with the wrapping implementation \cite{CandesFDCT05}. The Symmlet wavelet with 6 vanishing moments was used throughout all experiments. For each transform, two images were tested Barbara ($512 \times 512$) and Peppers ($256 \times 256$), and each image, was contaminated with zero-mean white Gaussian noise with increasing standard deviation $\sigma \in \{5,10,15,20,25,30\}$, corresponding to input PSNR values $\{34.15,28.13,24.61,22.11,20.17,18.59,14.15\}$ dB. At each combination of test image and noise level, ten noisy versions were generated. Then, block denoising was ten applied to each of the ten noisy images for each block size $L \in \{1,2,4,8,16\}$ and threshold $\lambda \in \{2,3,4,4.5,5,6\}$, and the average output PSNR over the ten realizations was computed. This yields one plot of average output PSNR as a function of $\lambda$ and $L$ at each combination (image-noise level-transform). The results are depicted in Fig.\ref{fig:dwtimpact}, Fig.\ref{fig:udwtimpact} and Fig.\ref{fig:curvimpact} for respectively the DWT, UDWT and FDCT. One can see that the maximum of PSNR occurs at $L=4$ (for $\lambda \geq 3$) whatever the transform and image, and this value turns to be the choice dictated by the theoretical procedure. As far as the influence of $\lambda$ is concerned, the PSNR attains its exact highest peak at different values of $\lambda$ depending on the image, transform and noise level. For the DWT, this maximum PSNR takes place near the theoretical threshold $\lambda_* \approx 4.5$ as expected from  Proposition \ref{mich2}. Even with the other redundant transforms, that correspond to tight frames for which Proposition \ref{mich2} is not rigorously valid, a sort of plateau is reached near $\lambda=4.5$. Only a minor improvement can be gained by taking a higher threshold $\lambda$; see e.g. Fig.\ref{fig:udwtimpact} or \ref{fig:curvimpact} with Peppers for $\sigma \geq 20$. Note that this improvement by taking a higher $\lambda$ for redundant transforms (i.e. non i.i.d. Gaussian noise) is formally predicted by Proposition \ref{mich}. Even though the estimate of Proposition \ref{mich} was expected to be rather crude. To summarize, the value $4.50524...$ intended to work for orthobases seems to yield good results also with redundant transforms.

\subsection{Comparative study}
\begin{figure}
\hspace*{-2cm}
\begin{tabular}{@{}c@{}c@{}}
Barbara $512 \times 512$ & Lena $512\times 512$ \\
\includegraphics[width=0.65\textwidth]{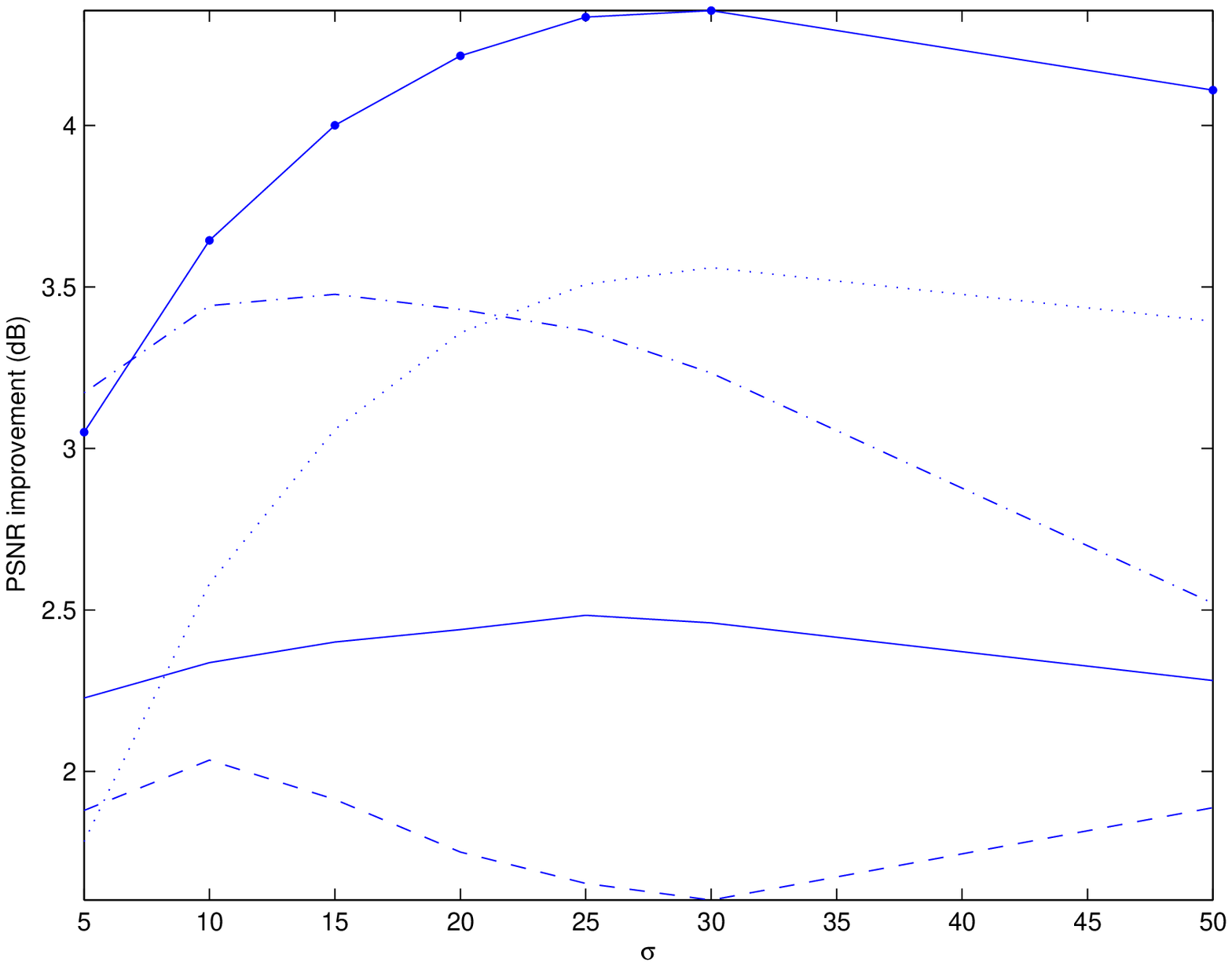}&\includegraphics[width=0.65\textwidth]{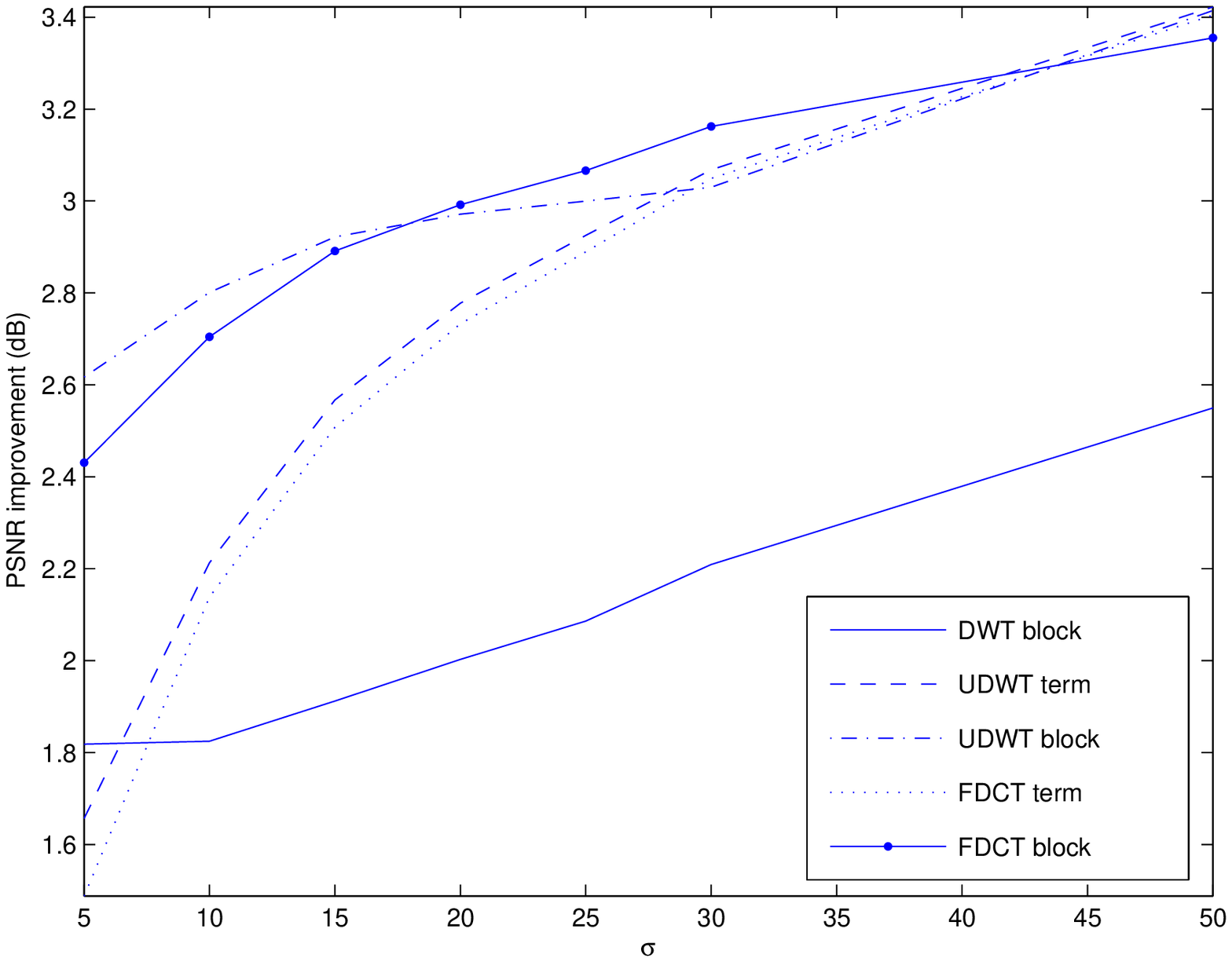} \\
House $256 \times 256$ & Boat $512\times 512$ \\
\includegraphics[width=0.65\textwidth]{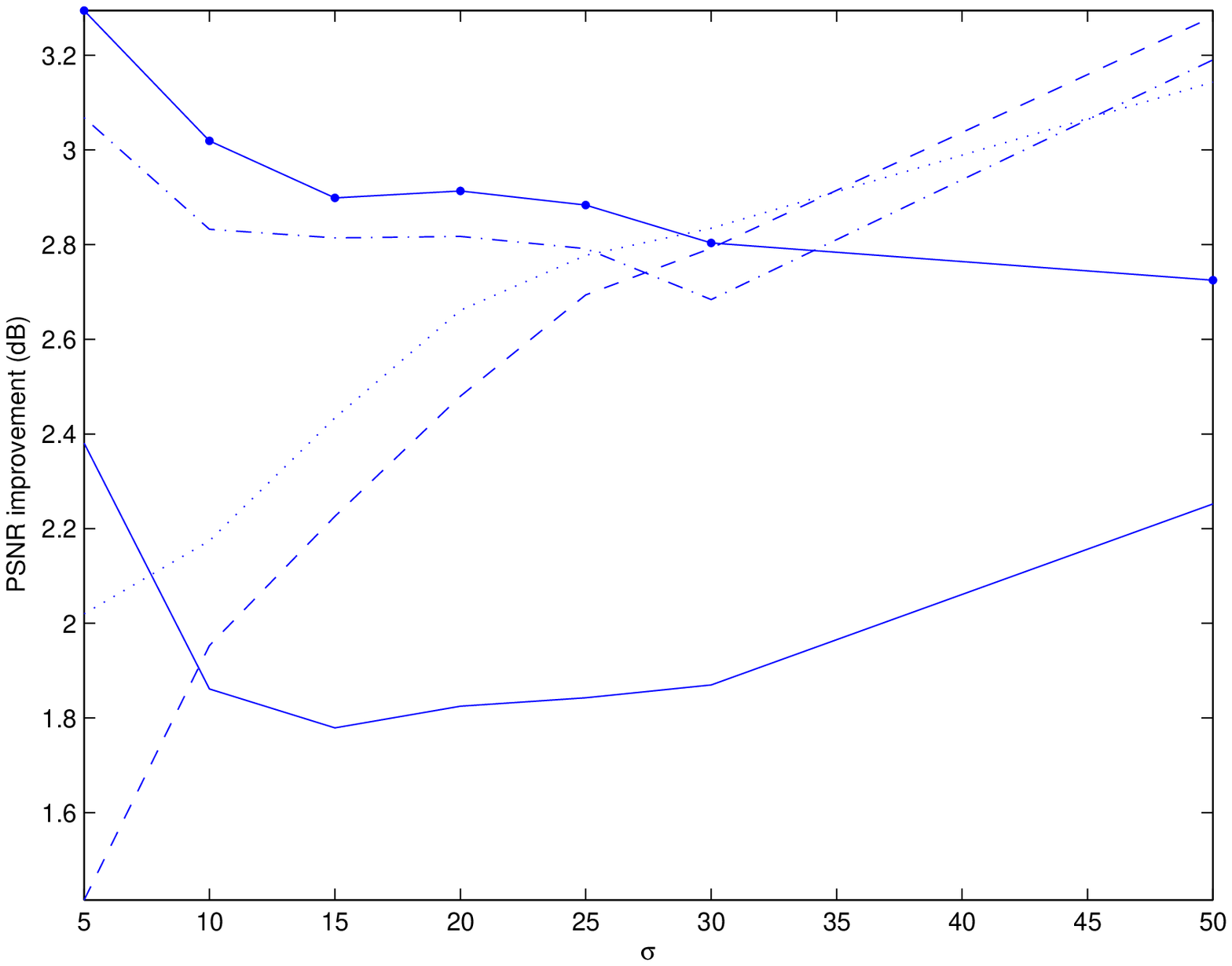}&\includegraphics[width=0.65\textwidth]{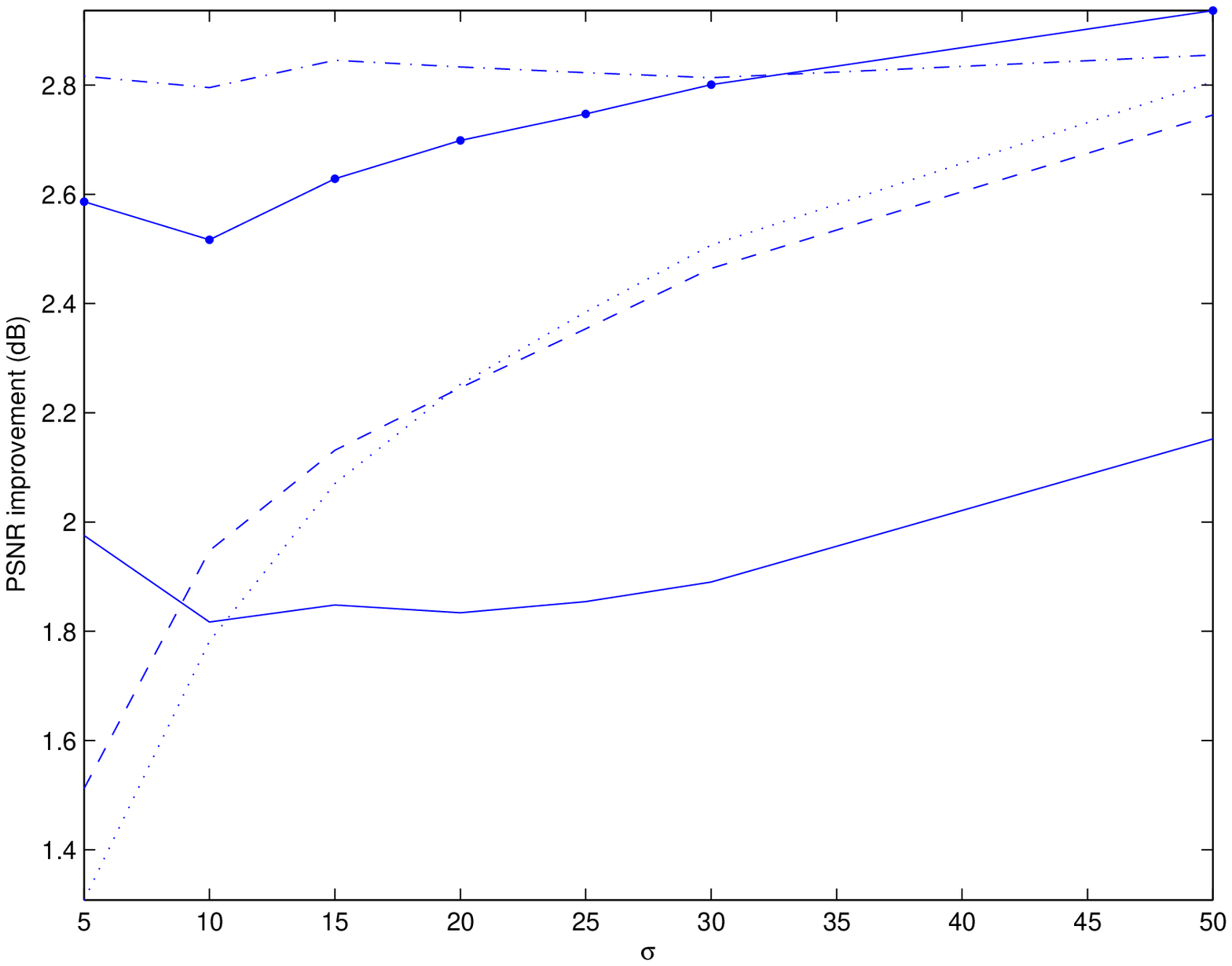} \\
Fingerprint $512\times 512$ & Peppers $256 \times 256$ \\
\includegraphics[width=0.65\textwidth]{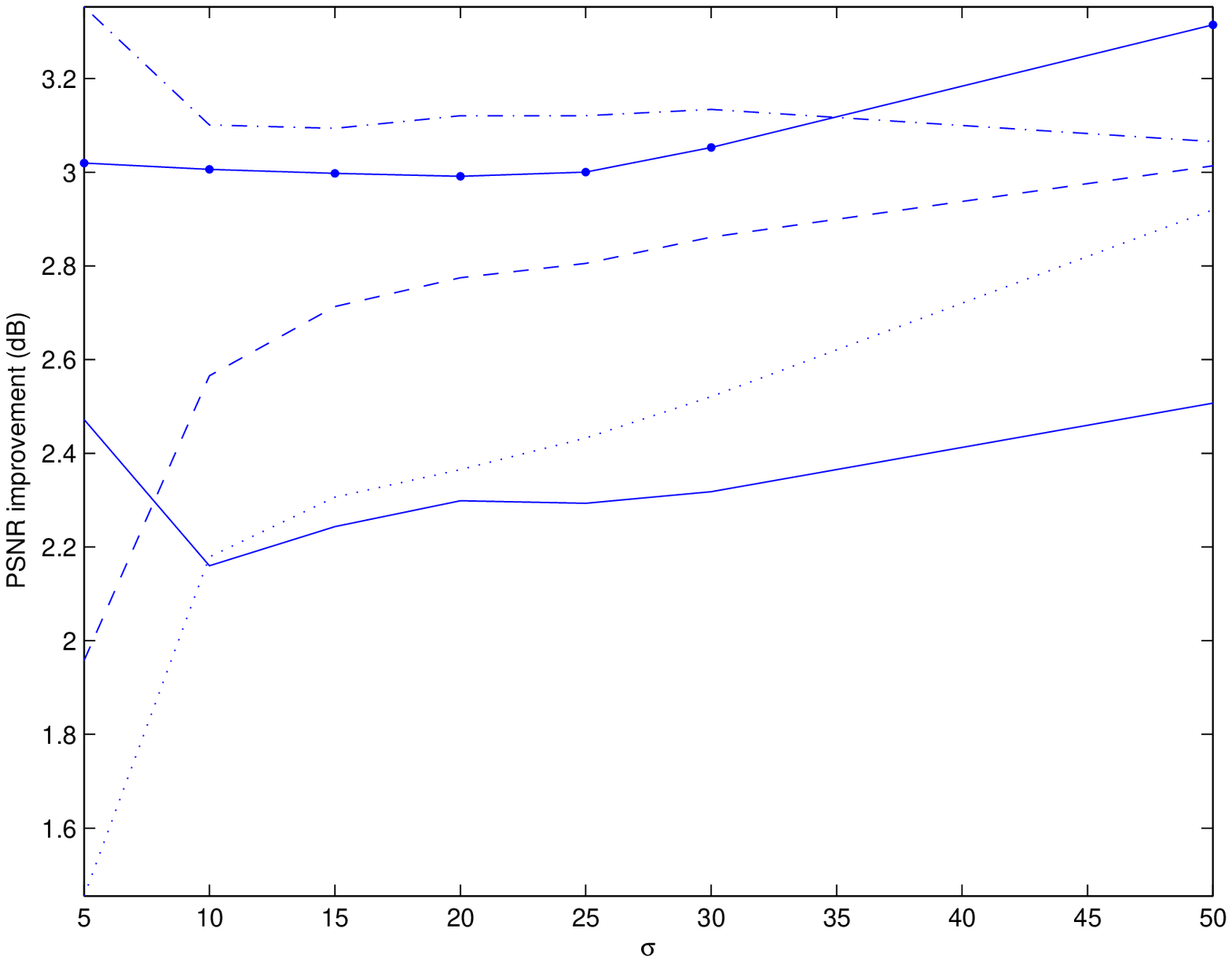}&\includegraphics[width=0.65\textwidth]{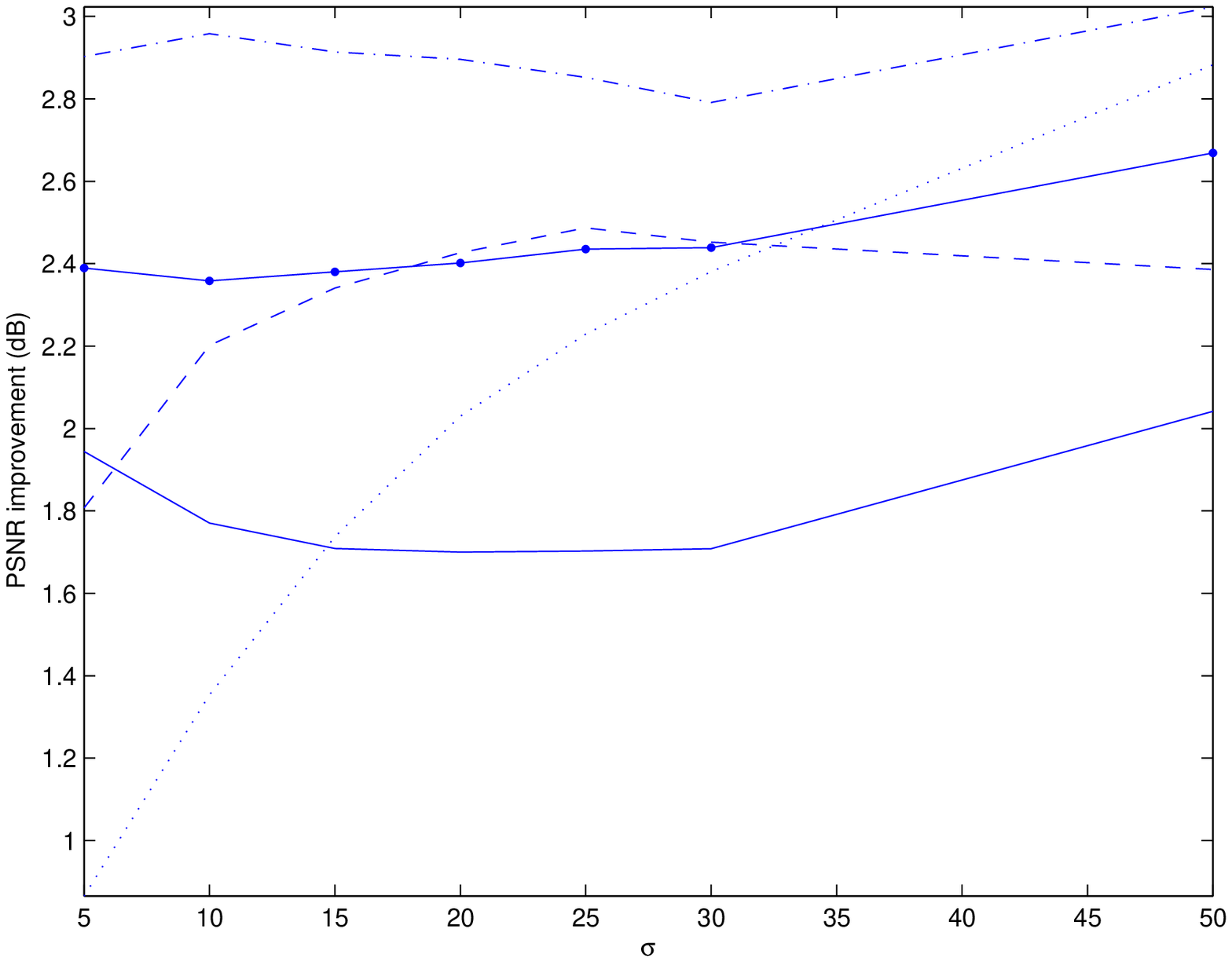}
\end{tabular}
\caption{Block vs term-by-term thresholding. Each plot corresponds to PSNR improvement over DWT term-by-term thresholding as a function of $\sigma$.}
\label{fig:blockvsterm}
\end{figure}

\paragraph*{Block vs term-by-term}
It is instructive to quantify the improvement brought by block denoising compared to term-by-term thresholding. For reliable comparison, we applied the denoising algorithms to six standard grayscale images with different contents of size $512 \times 512$ (Barbara, Lena, Boat and Fingerprint) and $256 \times 256$ (House and Peppers). All images were normalized to a maximum grayscale value 255. The images were corrupted by a zero-mean white Gaussian noise with standard deviation $\sigma \in \{5,10,15,20,25,30\}$. The output PSNR was averaged over ten realizations, and all algorithms were applied to the same noisy versions. The threshold used with individual thresholding was set to the classical value $3\sigma$ for the (orthogonal) DWT, and $3\sigma$ for all scales and $4\sigma$ at the finest scale for the (redundant) UDWT and FDCT. The results are displayed in Fig.\ref{fig:blockvsterm}. Each plot corresponds to PSNR improvement over DWT term-by-term thresholding as a function of $\sigma$. To summarize,
\begin{itemize}
\item Block shrinkage improves the denoising results in general compared to individual thresholding. Even though the improvement extent decreases with increasing $\sigma$. The PSNR increase brought by block denoising with a given transform compared to individual thresholding with the same transform can be up to 2.55 dB.
\item Owing to block shrinkage, even the orthogonal DWT becomes competitive with redundant transforms. For Barbara, block denoising with DWT is even better than individual thresholding in the translation-invariant UDWT.
\item For some images (e.g. Peppers or House), block denoising with curvelets can be slightly outperformed by its term-by-term thresholding counterpart for $\sigma=50$.
\item As expected, no transform is the best for all images. Block denoising with curvelets is more beneficial to images with high frequency content (e.g. anisotropic oscillating patterns in Barbara). For the other images, and except Peppers, block denoising with UDWT or curvelets are comparable ($\sim0.2$ dB difference).
\end{itemize}
Note that the additional computational burden of block shrinkage compared to individual thresholding is limited: respectively 0.1s, 1s and 0.7s for the DWT, UDWT and FDCT with $512 \times 512$ images, and less than 0.03s, 0.2s and 0.1 for $256 \times 256$ images. The algorithms were run under Matlab with an Intel Xeon 3GHz CPU, 8Gb RAM.

\begin{figure}
\centerline{\includegraphics{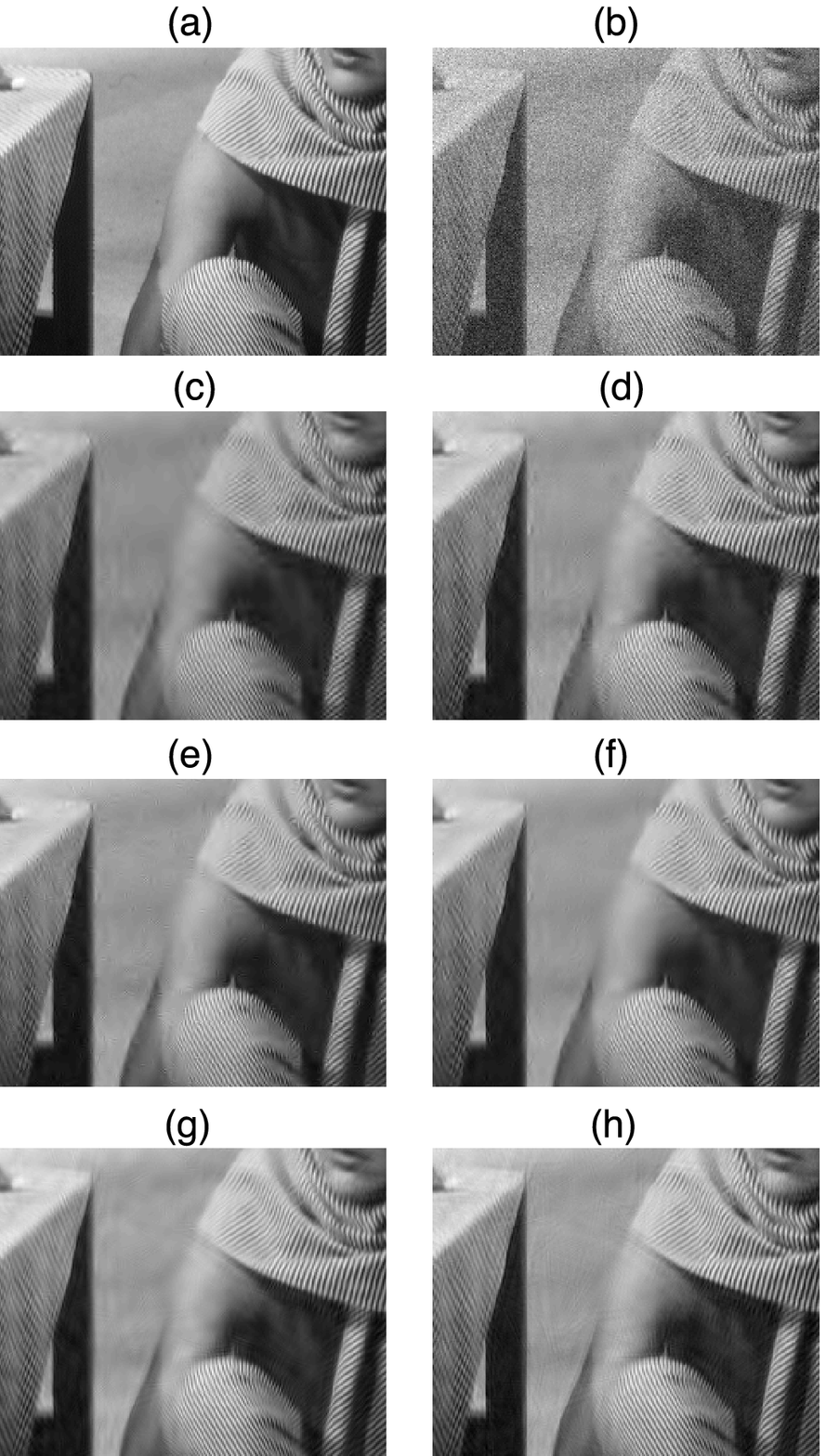}}
\caption{Visual comparison of our block denoising to BLS-GSM on Barbara $512 \times 512$. (a) original. (b) noisy $\sigma=20$. (c), (e) and (g) block denoising with respectively DWT (28.04 dB), UDWT (29.01 dB) and FDCT (30 dB). (d), (f) and (h) BLS-GSM with respectively DWT (28.6 dB), UDWT (29.3 dB) and FDCT (30.07 dB).}
\label{fig:barbara}
\end{figure}

\begin{figure}
\centerline{\includegraphics{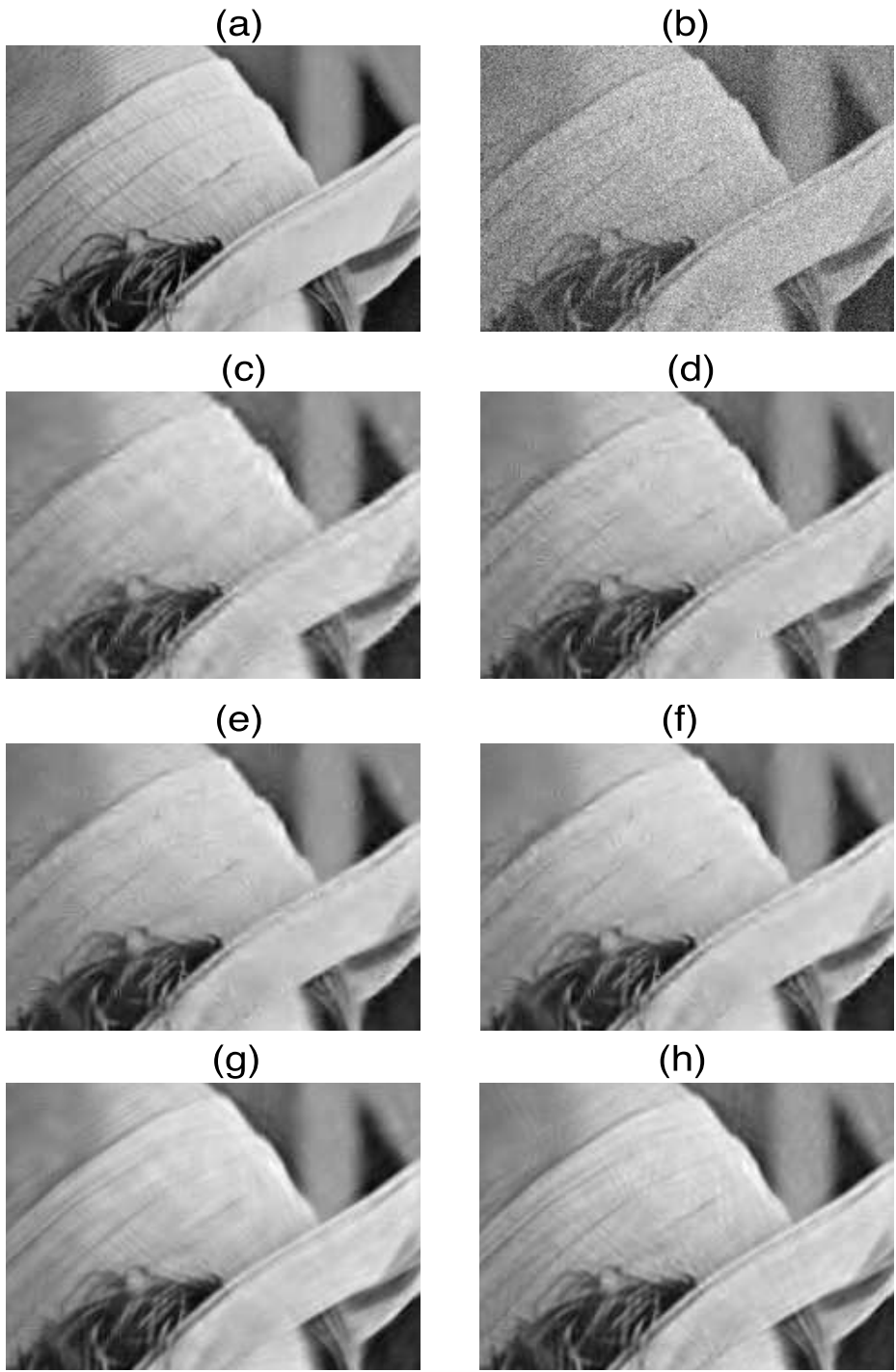}}
\caption{Visual comparison of our block denoising to BLS-GSM on Lena $512 \times 512$. (a) original. (b) noisy $\sigma=20$. (c), (e) and (g) block denoising with respectively DWT (30.51 dB), UDWT (31.47 dB) and FDCT (31.48 dB). (d), (f) and (h) BLS-GSM with respectively DWT (30.62 dB), UDWT (32 dB) and FDCT (31.6 dB).}
\label{fig:lena}
\end{figure}

\paragraph*{Block vs BLS-GSM}
The described block denoising procedure has been compared to one of state-of-the-art denoising methods in the literature BLS-GSM \cite{Portilla03}. BLS-GSM is a widely used reference in image denoising experiments reported in the literature. BLS-GSM uses a sophisticated prior model of the joint distribution within each block of coefficients, and then computes the Bayesian posterior conditional mean estimator by numerical integration. For fair comparison, BLS-GSM was also adapted and implemented with the curvelet transform. The two algorithms were applied to the same ten realizations of additive white Gaussian noise with $\sigma$ in the same range as before. The output PSNR values averaged over the ten realizations for each of the six tested image are tabulated in Table \ref{tab:compare}. By inspection of this table, the performance of block denoising and BLS-GSM remain comparable whatever the transform and image. None of them outperforms the other for all transforms and all images. When comparing both algorithms for the DWT transform, the maximum difference between the corresponding PSNR values is $0.5$ dB in favor of block shrinkage. For the UDWT and FDCT, the maximum difference is $\sim 0.6$ dB in BLS advantage. Visual inspection of Fig.\ref{fig:barbara} and \ref{fig:lena} is in agreement with the quantitative study we have just discussed. For each transform, differences between the two denoisers are hardly visible. Our procedure is however much simpler to implement and has a much lower computational cost than BLS-GSM as can be seen from Table~\ref{tab:time}. Our algorithm can be up to 10 times faster than BLS-GSM while reaching comparable denoising performance. As stated in the previous paragraph, the bulk of computation in our algorithm is essentially invested in computing the forward and inverse transforms.

\begin{table}
\begin{tabular}{cc}
$512 \times 512$ image & $256 \times 256$ image \\
\begin{tabular}{|l|c|c|c|}\hline
&         	DWT        &        UDWT    &    FDCT     \\\hline\hline
Block  &  	0.22       &        2.6     &    5.8       \\\hline
BLS-GSM    &  	3    	   &        26      &    30	      \\\hline
\end{tabular} &
\begin{tabular}{|l|c|c|c|}\hline
&         	DWT          &        UDWT       &    FDCT     \\\hline\hline
Block  &  	0.045        &        0.45       &    1.2           \\\hline
BLS-GSM    &  	1            &        5.5        &    6.6           \\\hline
\end{tabular}
\end{tabular}
\caption{Execution times in seconds for $512 \times 512$ images and $256 \times 256$ images. The algorithms were run under Matlab with an Intel Xeon 3GHz CPU, 8Gb RAM.}
\label{tab:time}
\end{table}

\begin{table}
\hspace*{-1cm}
\begin{scriptsize}
\hspace*{-1.1cm}
\begin{tabular}{@{}c@{}@{}c@{}}
Barbara $512 \times 512$ & Lena $512\times 512$ \\
\begin{scriptsize}
\begin{tabular}{|l|c|c|c|c|c|c|c|}
\hline
\textbf{$\sigma$}&5   &10   &15   &20   &25   &30   &50   \\
\textbf{\psnrin}&34.15&28.13&24.61&22.11&20.17&18.59&14.15\\\hline\hline
\textbf{Block DWT}&36.81&32.50&30.07&28.41&27.16&26.16&23.74\\\hline
\textbf{BLS-GSM DWT}&36.87&32.65&30.26&28.61&27.40&26.40&23.90\\\hline
\textbf{Block UDWT}&37.37&33.24&30.80&29.09&27.77&26.70&24.01\\\hline
\textbf{BLS-GSM UDWT}&37.44&33.43&31.06&29.40&28.16&27.13&24.49\\\hline
\textbf{Block FDCT}&37.57&33.68&31.52&30.00   &28.83&27.86&25.38\\\hline
\textbf{BLS-GSM FDCT}&37.63&33.82&31.64&30.08&28.93&28.01&25.36\\\hline
\end{tabular}
\end{scriptsize}&
\begin{scriptsize}
\begin{tabular}{|c|c|c|c|c|c|c|}
\hline
5   &10   &15   &20   &25   &30   &50   \\
34.15&28.13&24.61&22.11&20.17&18.59&14.15\\\hline\hline
37.61&34.05&31.99&30.62&29.58&28.71&26.36\\\hline
37.41&33.97&31.68&30.62&29.62&28.70&26.36\\\hline
38.02&34.75&32.85&31.48&30.41&29.53&27.16\\\hline
38.16&35.15&33.34&32.02&30.97&30.13&27.78\\\hline
38.09&34.78&32.86&31.45&30.43&29.55&27.12\\\hline
38.10&34.93&33.03&31.60&30.53&29.65&27.02\\\hline
\end{tabular}
\end{scriptsize}\\
House $256 \times 256$ & Boat $512\times 512$ \\
\begin{scriptsize}
\begin{tabular}{|l|c|c|c|c|c|c|c|}
\hline
\textbf{$\sigma$}&5   &10   &15   &20   &25   &30   &50   \\
\textbf{\psnrin}&34.15&28.13&24.61&22.11&20.17&18.59&14.15\\\hline\hline
\textbf{Block DWT}&37.63&33.47&31.33&29.86&28.76&27.79&25.41\\\hline
\textbf{BLS-GSM DWT}&37.43&33.97&31.77&29.88&29.17&28.43&26.12\\\hline
\textbf{Block UDWT}&38.10&34.31&32.31&30.86&29.75&28.80&26.35\\\hline
\textbf{BLS-GSM UDWT}&38.17&34.79&32.95&31.52&30.41&29.49&27.00\\\hline
\textbf{Block FDCT}&38.35&34.36&32.04&30.32&29.70&28.71&25.90\\\hline
\textbf{BLS-GSM FDCT}&38.47&34.69&32.47&30.92&29.71&28.72&25.93\\\hline
\end{tabular}
\end{scriptsize}&
\begin{scriptsize}
\begin{tabular}{|c|c|c|c|c|c|c|}
\hline
5   &10   &15   &20   &25   &30   &50   \\
34.15&28.13&24.61&22.11&20.17&18.59&14.15\\\hline\hline
36.41&32.52&30.41&28.93&27.81&26.97&24.83\\\hline
36.06&32.36&30.36&29.04&27.35&26.76&24.86\\\hline
36.89&33.15&31.11&29.67&28.59&27.71&25.45\\\hline
36.85&33.46&31.52&30.14&29.09&28.22&26.00\\\hline
36.89&33.07&31.03&29.65&28.59&27.70&25.49\\\hline
36.74&33.17&31.20&29.80&28.77&27.88&25.52\\\hline
\end{tabular}
\end{scriptsize}\\
Fingerprint $512\times 512$ & Peppers $256 \times 256$ \\
\begin{scriptsize}
\begin{tabular}{|l|c|c|c|c|c|c|c|}
\hline
\textbf{$\sigma$}&5   &10   &15   &20   &25   &30   &50   \\
\textbf{\psnrin}&34.15&28.13&24.61&22.11&20.17&18.59&14.15\\\hline\hline
\textbf{Block DWT}&35.74&31.37&29.10&27.53&26.33&25.34&22.84\\\hline
\textbf{BLS-GSM DWT}&35.53&31.08&28.82&27.08&26.01&25.11&22.72\\\hline
\textbf{Block UDWT}&36.22&31.89&29.62&28.06&26.87&25.90&23.37\\\hline
\textbf{BLS-GSM UDWT}&36.54&32.23&29.91&28.36&27.20&26.30&23.85\\\hline
\textbf{Block FDCT}&36.13&31.98&29.66&28.03&26.84&25.92&23.51\\\hline
\textbf{BLS-GSM FDCT}&36.34&32.14&29.82&28.21&27.05&26.14&23.70\\\hline
\end{tabular}
\end{scriptsize}&
\begin{scriptsize}
\begin{tabular}{|c|c|c|c|c|c|c|}
\hline
5   &10   &15   &20   &25   &30   &50   \\
34.15&28.13&24.61&22.11&20.17&18.59&14.15\\\hline\hline
36.81&32.56&30.28&28.64&27.42&26.42&23.77\\\hline
36.69&32.50&30.38&28.90&27.65&26.70&23.55\\\hline
37.48&33.60&31.37&29.74&28.52&27.52&24.71\\\hline
37.59&33.96&31.78&30.17&28.99&27.97&25.16\\\hline
37.09&33.14&30.86&29.17&28.01&27.09&24.38\\\hline
37.15&33.32&31.10&29.44&28.19&26.85&24.27\\\hline
\end{tabular}
\end{scriptsize}\end{tabular}
\end{scriptsize}

\caption{Comparison of average PSNR over ten realizations of block denoising and BLS-GSM, with three transforms.}
\label{tab:compare}
\end{table}

\subsection{Reproducible research}
Following the philosophy of reproducible research, a toolbox is made available freely for download at the address

{\centerline{\texttt{http://www.greyc.ensicaen.fr/$\sim$jfadili/software.html}}}

This toolbox is a collection of Matlab functions, scripts and datasets for image block denoising. It requires at least WaveLab 8.02 \cite{BuckheitDonoho95} to run properly. The toolbox implements the proposed block denoising procedure with several transforms and contains all scripts to reproduce the figures and tables reported in this paper.

\section{Conclusion}\label{conclusion}
In this paper, an Stein block thresholding algorithm for denoising $d$-dimensional data is proposed with a particular focus on 2D image. Our block denoising is a generalization of one-dimensional BlockJS to $d$ dimensions, with other transforms that orthogonal wavelets, and handles noise in the coefficient domain beyond the i.i.d. Gaussian case. Its minimax properties are investigated, and a fast and appealing algorithm is described. The practical performance of the designed denoiser were shown to be very promising with several transforms and a variety of test images. It turns out that the proposed block denoiser is much faster than state-of-the art competitors in the literature while reaching comparable denoising performance.

We believe however that there is still room for improvement of our procedure. For instance, for $d=2$, it would be interesting to investigate both theoretically and in practice how our results can be adapted to anisotropic blocks with possibly varying sizes. The rationale behind such a modification is to adapt the blocks to the geometry of the neighborhood. We expect that the analysis in this case, if possible, would be much more involved. Another interesting line of research would be to try to improve our convergence rates by relaxing the condition $q \leq p \wedge 2$. At this moment, given our definition of the smoothness space and our derivations in the proof (see appendix), we have not found a way around it yet. As remarked in subsection \ref{subsubsec:noiseseq}, a parameter $\delta$ was introduced, whose role becomes of interest when addressing linear inverse problems such as deconvolution. Extension of BlockJS to linear inverse problems remains also an open question. All these aspects need further investigation that we leave for a future work.


\appendix
\section*{Appendix: Proofs}
\label{proofs}
In this section, $C$ represents a positive constant which may differ from one term to another. We suppose that $n$ is large enough.

\section{Proof of Theorem \ref{mini2}}
We have the decomposition:
\begin{equation}\label{risk}
R(\widehat \theta^*,\theta) = R_1+R_2+R_3,
\end{equation}
where
$$R_1= \sum_{j=0}^{j_0-1}\sum_{\ell\in B_j}\sum_{{\bf k} \in D_j}\mathbb{E}\left((\widehat \theta^*_{j,\ell,{\bf k}}-\theta_{j,\ell,{\bf k}})^2\right), \ \ \ \ \ R_2=\sum_{j=j_0}^{J_*} \sum_{\ell\in B_j}\sum_{{\bf k} \in D_j }\mathbb{E}\left((\widehat \theta^*_{j,\ell,{\bf k}}-\theta_{j,\ell,{\bf k}})^2\right), $$
$$R_3=\sum_{j=J_*+1}^{\infty}\sum_{\ell\in B_j}\sum_{{\bf k} \in D_j} \theta_{j,\ell,{\bf k}}^2.$$
Let us bound the terms $R_1$, $R_3$ and $R_2$ (by order of difficulty).

{\bf The upper bound for $R_1$.} It follows from $(\mathrm{A}1)$ that
\begin{eqnarray}\label{a1}
R_1& = & n^{-r}\sum_{j=0}^{j_0-1} \sum_{\ell\in B_j}  \sum_{{\bf k} \in D_j}\mathbb{E}\left( z^2_{j,\ell,{\bf k}}\right)\le  Q_1 n^{-r}\sum_{j=0}^{j_0-1}2^{j(d_*+\delta)}\card(B_j) \nonumber \\
& =
& c_* Q_1 n^{-r}\sum_{j=0}^{j_0-1}2^{j(d_*+\delta+\upsilon) }  \le  C 2^{j_0(d_*+\delta+\upsilon)}n^{-r} \nonumber \\
& \le
&    C L^{(1/\min_{i=1,...,d} \mu_i))(d_*+\delta+\upsilon)}n^{-r}\le C (\log n)^{(1/(d\min_{i=1,...,d} \mu_i))(d_*+\delta+\upsilon)}n^{-r}  \nonumber \\
& \le
&  C n^{-2sr/(2s+\delta+d_*+\upsilon)}.
\end{eqnarray}
We used the inequality $2s/(2s+\delta+d_*+\upsilon)<1$ which implies that, for a large enough $n$, $(\log n)^{(1/(d\min_{i=1,...,d} \mu_i))(d_*+\delta+\upsilon)}n^{(2s/(2s+\delta+d_*+\upsilon)-1)r}\le 1$. \\

{\bf The upper bound for $R_3$.} We distinguish the case $q\le 2 \le p$ and the case $q\le p<2$.

For $q\le 2 \le p$, we have ${\boldsymbol\Theta}_{p,q}^s(M) \subseteq {\boldsymbol\Theta}_{2,q}^s(M)\subseteq {\boldsymbol\Theta}_{2,2}^s(M) $. Hence
\begin{eqnarray}\label{a31}
R_3 \le  M^2 \sum_{j=J_*+1}^{\infty}2^{-2js}\le C 2^{-2J_*s}\le C n^{-2sr/(d_*+\delta+\upsilon)} \le  C n^{-2sr/(2s+\delta+d_*+\upsilon)}. \nonumber\\
&
&
\end{eqnarray}
For $q\le p< 2$, we have ${\boldsymbol\Theta}_{p,q}^s(M) \subseteq {\boldsymbol\Theta}_{2,q}^{s-d_*/p+d_*/2}(M)  \subseteq {\boldsymbol\Theta}_{2,2}^{s-d_*/p+d_*/2}(M)$. We have
\begin{alignat*}{3}
&\quad \quad &s/(2s+\delta+d_*+\upsilon) &\le (s-d_*/p+d_*/2)/(d_*+\delta+\upsilon) \\
&\Leftrightarrow \quad &s(d_*+\delta+\upsilon) &\le (s-d_*/p+d_*/2)(2s+\delta+d_*+\upsilon) \\
&\Leftrightarrow \quad &0 &\le 2s^2-(d_*/p-d_*/2)(2s+\delta+d_*+\upsilon) \\
&\Leftrightarrow \quad &0 &\le 2s(s-d_*/p)+sd_*-(d_*/p-d_*/2)(\delta+d_*+\upsilon) ~.
\end{alignat*}
This implies that, if $sp>d_*$ and $s>(1/p-1/2)(\delta+d_*+\upsilon)$, we have $s/(2s+\delta+d_*+\upsilon)\le (s-d_*/p+d_*/2)/(d_*+\delta+\upsilon)$. Therefore,
\begin{eqnarray}\label{a32}
R_3& \le  & M^2 \sum_{j=J_*+1}^{\infty} 2^{-2j(s-d_*/p+d_*/2) }   \le C 2^{-2J_*(s-d_*/p+d_*/2)}\nonumber \\
& \le
&  C n^{-2r(s-d_*/p+d_*/2)/(d_*+\delta+\upsilon)} \le  C n^{-2sr/(2s+\delta+d_*+\upsilon)}.
\end{eqnarray}
Putting \eqref{a31} and \eqref{a32} together, we obtain the desired upper bound.

{\bf The upper bound for $R_2$.} We need the following result which will be proved later.
\begin{lemma} \label{mal}
Let $(v_i)_{i\in \mathbb{N}^*}$ be a sequence of real numbers and $(w_i)_{i\in \mathbb{N}^*}$ be a sequence of random variables. Set, for any $i\in\mathbb{N}^*$,
$$u_i=v_i+w_i.$$
Then, for any $m\in\mathbb{N}^*$ and any $\lambda>0$, the sequence of estimates $(\tilde{u}_i)_{i=1,...,m}$ defined by $\tilde{u}_i=u_i\left( 1-\lambda^2 \left( \sum_{i=1}^m u^2_i  \right)^{-1}\right)_+$
satisfies
$$\sum_{i=1}^m (\tilde{u}_i-v_i)^2 \le 10 \sum_{i=1}^m w_i^2\indic{\left\lbrace  \left(\sum_{i=1}^m w_i^2\right)^{1/2} > \lambda/2 \right\rbrace} + 10\min \left(  \sum_{i=1}^m v_i^2 , \lambda^2/4 \right) .$$
\end{lemma}
Lemma \ref{mal} yields
\begin{equation}\label{a2}
R_2 =  \sum_{j=j_0}^{J_*}\sum_{\ell\in B_j}\sum_{{\bf K} \in \mathcal{A}_j}\sum_{{\bf k} \in U_{j,{\bf K}}}\mathbb{E}\left((\widehat \theta^*_{j,\ell,{\bf k}}-\theta_{j,\ell,{\bf k}})^2\right) \le 10 (B_1+B_2),
\end{equation}
where
$$B_1=n^{-r}\sum_{j=j_0}^{J_*}\sum_{\ell\in B_j}\sum_{{\bf K} \in \mathcal{A}_j}  \sum_{{\bf k} \in U_{j,{\bf K}}}\mathbb{E}\left( z_{j,\ell,{\bf k}} ^2\indic{\left\lbrace \sum_{{\bf k} \in U_{j,{\bf K}}} z_{j,\ell,{\bf k}}^2 > \lambda_*2^{\delta j}L^d/4 \right\rbrace}\right)$$
and
$$B_2=\sum_{j=j_0}^{J_*}\sum_{\ell\in B_j}\sum_{{\bf K} \in \mathcal{A}_j}\min \left(  \sum_{{\bf k} \in U_{j,{\bf K}}} \theta_{j,\ell,{\bf k}}^2 , \lambda_* 2^{\delta j}L^d n^{-r}/4 \right) .$$
Using $(\mathrm{A}2)$, we bound $B_1$ by
\begin{equation}\label{b1}
B_1\le Q_2 n^{-r} \le Q_2   n^{-2sr/(2s+\delta +d_*+\upsilon)}.
\end{equation}
To bound $B_2$, we again distinguish the case $q\le 2 \le p$ and the case $q\le p<2$.

For $q\le 2 \le p$, we have $\theta \in {\boldsymbol\Theta}_{p,q}^s(M) \subseteq {\boldsymbol\Theta}_{2,q}^s(M) \subseteq {\boldsymbol\Theta}_{2,2}^s(M) $.
Let $j_s$ be the integer $j_s=\lfloor (r/(2s+\delta+d_*+\upsilon)) \log_2 n\rfloor $. We then obtain the bound
\begin{eqnarray}\label{b2}
B_2& \le & 4^{-1} \lambda_*L^d n^{-r} \sum_{j=j_0}^{j_s} 2^{j\delta}\card (\mathcal{A}_j)\card(B_j)+\sum_{j=j_s+1}^{J_*}\sum_{\ell\in B_j} \sum_{{\bf k} \in D_j} \theta_{j,\ell,{\bf k}}^2\nonumber \\
& \le
& 4^{-1}c_* \lambda_*L^d n^{-r}  \sum_{j=j_0}^{j_s} 2^{j (d_*+\delta+\upsilon)} L^{-d}+\sum_{j=j_s+1}^{J_*}\sum_{\ell\in B_j} \sum_{{\bf k} \in D_j} \theta_{j,\ell,{\bf k}}^2\nonumber \\
& \le
& C n^{-r} 2^{j_s (d_*+\delta+\upsilon)} +M^2 \sum_{j=j_s+1}^{J_*} 2^{-2js}\nonumber \\
& \le
& C n^{-r} 2^{j_s (d_*+\delta+\upsilon)} + C 2^{-2j_ss} \le C n^{-2sr/(2s+\delta+d_*+\upsilon)}.
\end{eqnarray}
Putting \eqref{a2}, \eqref{b1} and \eqref{b2} together, it follows immediately that
\begin{equation}\label{a2fin}
R_2\le C n^{-2sr/(2s+\delta+d_*+\upsilon)}.
\end{equation}
Let's now turn to the case $q\le  p<2$. Let $j^*_s$ be the integer $j^*_s=\lfloor (r/(2s+\delta+d_*+\upsilon)) \log_2 (n/\log n)\rfloor$. We have
\begin{equation}\label{b22}
B_2\le D_{1}+D_{2}+D_{3},
\end{equation}
where
$$D_1=4^{-1} \lambda_*L^d n^{-r} \sum_{j=j_0}^{j^*_s} 2^{j\delta}\card(\mathcal{A}_j)\card(B_j),$$
$$D_2=4^{-1} \lambda_*L^d n^{-r} \sum_{j=j^*_s+1}^{J_*}\sum_{\ell\in B_j} \sum_{{\bf K} \in \mathcal{A}_j} 2^{\delta j}\indic{\left\lbrace \sum_{{\bf k} \in U_{j,{\bf K}}} \theta_{j,\ell,{\bf k}}^2> \lambda_*2^{\delta j} L^d n^{-r}/4 \right\rbrace}$$
and
$$D_3= \sum_{j=j^*_s+1}^{J_*}\sum_{\ell\in B_j} \sum_{{\bf K} \in \mathcal{A}_j} \sum_{{\bf k} \in U_{j,{\bf K}}} \theta_{j,\ell,{\bf k}}^2\indic{\left\lbrace \sum_{{\bf k} \in U_{j,{\bf K}}} \theta_{j,\ell,{\bf k}}^2 \le \lambda_*2^{\delta j}L^d n^{-r}/4 \right\rbrace}.$$
We have
\begin{eqnarray}\label{d1}
D_1 &\le  &4^{-1}c_* \lambda_*L^d n^{-r}  \sum_{j=j_0}^{j^*_s} 2^{j (d_*+\delta+\upsilon)} L^{-d} \le  C n^{-r} 2^{j^*_s (d_*+\delta+\upsilon)}\nonumber \\
& \le
&  C (\log n/n)^{2sr/(2s+\delta+d_*+\upsilon)}.
\end{eqnarray}
Moreover, using the classical inequality $\|\theta\|_2^p \leq \|\theta\|_p^p, p < 2$, we obtain
\begin{eqnarray}\label{d2a}
D_2 &\le & C L^d n^{-r} (L^d n^{-r})^{-p/2} \sum_{j=j^*_s+1}^{J_*} 2^{\delta j(1-p/2)} \sum_{\ell\in B_j}\sum_{{\bf K} \in \mathcal{A}_j} \left(\sum_{{\bf k} \in U_{j,{\bf K}}} \theta_{j,\ell,{\bf k}}^2 \right)^{p/2} \nonumber \\
& \le
& C(\log n/n)^{r(1-p/2)}\sum_{j=j^*_s+1}^{J_*} 2^{\delta j(1-p/2)} \sum_{\ell\in B_j}\sum_{{\bf k} \in D_j}|\theta_{j,\ell,{\bf k}}|^p.
\end{eqnarray}
Since $q \le  p$, we have ${\boldsymbol\Theta}_{p,q}^s(M) \subseteq {\boldsymbol\Theta}_{p,p}^{s}(M)$. Combining this with $sp>d_*$ and $s>(1/p-1/2)(\delta+d_*+\upsilon)$, we obtain
\begin{eqnarray}\label{d2}
D_2 & \le  &  C(\log n/n)^{r(1-p/2)}\sum_{j=j^*_s+1}^{J_*} 2^{\delta j(1-p/2)}2^{-j(s+d_*/2-d_*/p)p} \nonumber \\
& \le
& C(\log n/n)^{r(1-p/2)}2^{-j_s^*(s+d_*/2-d_*/p-\delta/p+\delta/2)p}\nonumber\\
& \le
& C (\log n/n)^{(2s+\upsilon(1-p/2))r/(2s+\delta +d_*+\upsilon)}\nonumber \\
& \le
&  C (\log n/n)^{2sr/(2s+\delta +d_*+\upsilon)}.
\end{eqnarray}

We have, for any ${\bf k}\in U_{j, {\bf K}}$, the inclusion $\left\lbrace \sum_{{\bf k} \in U_{j, {\bf K}}} \theta_{j,\ell,{\bf k}}^2 \le \lambda_*2^{\delta j}L^d n^{-r}/4 \right\rbrace\subseteq \left\lbrace |\theta_{j,\ell,{\bf k}}|\le (\lambda_*2^{\delta j}L^d n^{-r})^{1/2}/2 \right\rbrace$.
Therefore,
\begin{eqnarray*}
D_3& \le & \sum_{j=j^*_s+1}^{J_*}\sum_{\ell\in B_j} \sum_{{\bf K} \in \mathcal{A}_j} \sum_{{\bf k} \in U_{j,{\bf K}}} \theta_{j,\ell,{\bf k}}^2\indic{\left\lbrace |\theta_{j,\ell,{\bf k}}|\le (\lambda_*2^{\delta j}L^d n^{-r})^{1/2}/2 \right\rbrace} \\
&  \le
&  C (\lambda_*L^d n^{-r})^{1-p/2}\sum_{j=j^*_s+1}^{J_*}2^{j\delta(1-p/2)} \sum_{\ell\in B_j}\sum_{{\bf k} \in D_j} |\theta_{j,\ell,{\bf k}}|^p ~,
\end{eqnarray*}
which is the same bound as for $D_2$ in \eqref{d2a}. Then using similar arguments as those used for in \eqref{d2}, we arrive at
\begin{equation}\label{d3}
D_3 \le C (\log n/n)^{2sr/(2s+\delta +d_*+\upsilon)}.
\end{equation}
Inserting \eqref{d1}, \eqref{d2} and \eqref{d3} into \eqref{b22}, it follows that
\begin{equation}\label{a2fina}
R_2\le C (\log n/n)^{2sr/(2s+\delta +d_*+\upsilon)}.
\end{equation}
Finally, bringing \eqref{risk}, \eqref{a1}, \eqref{a31},\eqref{a32}, \eqref{a2fin} and \eqref{a2fina} together we obtain
$$\sup_{\theta \in {{\boldsymbol\Theta}}_{p,q}^s (M)} R(\widehat \theta^*,\theta)\le R_1+R_2+R_3\le C \rho_n,$$
where $\rho_n$ is defined by \eqref{rates}.
This ends the proof of Theorem \ref{mini2}.

\section{Proof of Lemma \ref{mal}}
We have
\begin{equation}\label{tru}
\sum_{i=1}^m (\tilde{u}_i-v_i)^2=\max \left(A,B\right),
\end{equation}
where
$$A=\sum_{i=1}^m \left(w_i -\lambda^2 u_i \left( \sum_{i=1}^m u^2_i \right)^{-1} \right)^2\indic{\left\lbrace \left( \sum_{i=1}^m u_i ^2 \right)^{1/2} > \lambda \right\rbrace}, \ \ \ B= \sum_{i=1}^m v_i^2\indic{\left\lbrace \left( \sum_{i=1}^m u_i^2\right)^{1/2} \le \lambda \right\rbrace}.$$
Let's bound $A$ and $B$, in turn.

{\bf The upper bound for $A$.} Using the elementary inequality $(a-b)^2 \le 2(a^2+b^2)$, we have
\begin{eqnarray}\label{ds}
A& \le & 2 \sum_{i=1}^m \left( w_i^2 +\lambda^4 u_i^2 \left( \sum_{i=1}^m u^2_i \right)^{-2} \right)\indic{\left\lbrace \left( \sum_{i=1}^m u_i ^2 \right)^{1/2} > \lambda \right\rbrace}\nonumber \\
& =
& 2 \left(\sum_{i=1}^m  w_i^2 +\lambda^4 \left( \sum_{i=1}^m u^2_i \right)^{-1} \right)\indic{\left\lbrace \left( \sum_{i=1}^m u_i ^2 \right)^{1/2} > \lambda \right\rbrace}\nonumber \\
& \le
& 2 \left(\sum_{i=1}^m  w_i^2 +\lambda^2 \right)\indic{\left\lbrace \left( \sum_{i=1}^m u_i ^2 \right)^{1/2} > \lambda \right\rbrace}.
\end{eqnarray}
Set $$D=2 \left(\sum_{i=1}^m  w_i^2 +\lambda^2 \right) \indic{\left\lbrace \left( \sum_{i=1}^m u_i ^2 \right)^{1/2} > \lambda \right\rbrace}.$$

We have the decomposition
\begin{equation}\label{ds1}
D= D_1+D_2,
\end{equation}
where
$$D_1=D \indic{\left\lbrace  \left(\sum_{i=1}^m w_i^2\right)^{1/2} > \lambda/2 \right\rbrace}, \ \ \ \ \ \ \ \ D_2=D\indic{\left\lbrace  \left(\sum_{i=1}^m w_i^2\right)^{1/2} \le \lambda/2 \right\rbrace}.$$
We clearly have
\begin{eqnarray}\label{ds2}
D_1  &\le &2 \left(\sum_{i=1}^m  w_i^2 +\lambda^2 \right) \indic{\left\lbrace  \left(\sum_{i=1}^m w_i^2\right)^{1/2} > \lambda/2 \right\rbrace} \le  10 \sum_{i=1}^m  w_i^2 \indic{\left\lbrace  \left(\sum_{i=1}^m w_i^2\right)^{1/2} > \lambda/2 \right\rbrace}.\nonumber \\
&
&
\end{eqnarray}
Using the Minkowski inequality, we have the inclusion $\left\lbrace  \left(\sum_{i=1}^m u_i^2\right)^{1/2} > \lambda \right\rbrace \cap \left\lbrace  \left(\sum_{i=1}^m w_i^2\right)^{1/2} \le  \lambda/2 \right\rbrace
 \subseteq
 \left\lbrace  \left(\sum_{i=1}^m v_i^2\right)^{1/2} > \lambda/2 \right \rbrace \cap \left\lbrace  \left(\sum_{i=1}^m w_i^2\right)^{1/2} \le  \lambda/2 \right\rbrace.$
Therefore
\begin{eqnarray}\label{ds3}
D_2 & \le
&  2 \left(\sum_{i=1}^m  w_i^2 +\lambda^2 \right) \indic{\left\lbrace  \left(\sum_{i=1}^m v_i^2\right)^{1/2} > \lambda/2 \right \rbrace \cap \left\lbrace  \left(\sum_{i=1}^m w_i^2\right)^{1/2} \le  \lambda/2 \right\rbrace  }\nonumber \\
&  \le
&   10 \min \left(  \sum_{i=1}^m v_i^2 , \lambda^2/4 \right).
\end{eqnarray}
If we combine \eqref{ds}, \eqref{ds1}, \eqref{ds2} and \eqref{ds3}, we obtain
\begin{equation}\label{aaa}
A\le D \le 10 \sum_{i=1}^m w_i^2\indic{\left\lbrace  \left(\sum_{i=1}^m w_i^2\right)^{1/2} > \lambda/2 \right\rbrace}+ 10 \min \left(  \sum_{i=1}^m v_i^2 , \lambda^2/4 \right).
\end{equation}
{\bf The upper bound for $B$.} We have the decomposition
\begin{equation}\label{bb}
B= G_1+G_2
\end{equation}
$$G_1=B \indic{\left\lbrace  \left(\sum_{i=1}^m w_i^2\right)^{1/2} > \lambda/2 \right\rbrace}, \ \ \ \ \ \ G_2=B\indic{\left\lbrace  \left(\sum_{i=1}^m w_i^2\right)^{1/2} \le \lambda/2 \right\rbrace}.$$
Using the Minkowski inequality, we have again the inclusion $\left\lbrace  \left(\sum_{i=1}^m u_i^2\right)^{1/2} \le \lambda \right\rbrace \cap \left\lbrace  \left(\sum_{i=1}^m w_i^2\right)^{1/2} >  \lambda/2 \right\rbrace
 \subseteq
 \left\lbrace  \left(\sum_{i=1}^m v_i^2\right)^{1/2} \le 3 \left(\sum_{i=1}^m w_i^2\right)^{1/2} \right \rbrace \cap \left\lbrace  \left(\sum_{i=1}^m w_i^2\right)^{1/2} >  \lambda/2 \right\rbrace$.
It follows that
\begin{eqnarray}\label{g1}
G_1 & \le &  \sum_{i=1}^m v_i^2\indic{ \left\lbrace  \left(\sum_{i=1}^m v_i^2\right)^{1/2} \le 3 \left(\sum_{i=1}^m w_i^2\right)^{1/2} \right \rbrace \cap \left\lbrace  \left(\sum_{i=1}^m w_i^2\right)^{1/2} >  \lambda/2 \right\rbrace}\nonumber  \\
& \le
& 9  \sum_{i=1}^m w_i^2\indic{\left\lbrace  \left(\sum_{i=1}^m w_i^2\right)^{1/2} >  \lambda/2 \right\rbrace}.
\end{eqnarray}
Another application of the Minkowski inequality leads to the inclusion $\left\lbrace  \left(\sum_{i=1}^m u_i^2\right)^{1/2} \le \lambda \right\rbrace \cap \left\lbrace  \left(\sum_{i=1}^m w_i^2\right)^{1/2} \le  \lambda/2 \right\rbrace
\subseteq \left\lbrace  \left(\sum_{i=1}^m v_i^2\right)^{1/2} \le 3 \lambda/2  \right \rbrace \cap \left\lbrace  \left(\sum_{i=1}^m w_i^2\right)^{1/2} \le  \lambda/2 \right\rbrace$.
It follows that
\begin{eqnarray}\label{g2}
G_2 & \le &  \sum_{i=1}^m v_i^2\indic{\left\lbrace  \left(\sum_{i=1}^m v_i^2\right)^{1/2} \le 3 \lambda/2  \right \rbrace \cap \left\lbrace  \left(\sum_{i=1}^m w_i^2\right)^{1/2} \le  \lambda/2 \right\rbrace} \nonumber \\
& \le
& \min \left(  \sum_{i=1}^m v_i^2 , 9 \lambda^2/4 \right).
\end{eqnarray}
Therefore, if we combine \eqref{bb}, \eqref{g1} and \eqref{g2}, we obtain
\begin{equation}\label{bbb}
B \le 9 \sum_{i=1}^m w_i^2\indic{\left\lbrace  \left(\sum_{i=1}^m w_i^2\right)^{1/2} > \lambda/2 \right\rbrace}+\min \left(  \sum_{i=1}^m v_i^2 , 9 \lambda^2/4 \right).
\end{equation}
Putting \eqref{tru}, \eqref{aaa} and \eqref{bbb} together, we have
\begin{eqnarray*}
\sum_{i=1}^m (\tilde{u}_i-v_i)^2 & = & \max \left( A,B \right) \\
& \le
&  10 \sum_{i=1}^m w_i^2\indic{\left\lbrace  \left(\sum_{i=1}^m w_i^2\right)^{1/2} > \lambda/2 \right\rbrace} + 10\min \left(  \sum_{i=1}^m v_i^2 , \lambda^2/4 \right).
\end{eqnarray*}
Lemma \ref{mal} is proved.

\section{Proof of Proposition \ref{mich}}
First of all, notice that the Jensen inequality, $(\mathrm{A}3)$ and the fact that $\card(D_j) \leq 2^{jd_*}$ imply
\begin{eqnarray*}
\sup_{j\in \{0,...,J\}}\sup_{\ell\in B_j}2^{-j(d_*+\delta)} \sum_{{\bf k} \in D_j}\mathbb{E}\left( z^2_{j,\ell,{\bf k}}\right) & \le & \sup_{j\in \{0,...,J\}}2^{-j(d_*+\delta)}\sup_{\ell\in B_j} \sum_{{\bf k} \in D_j} \left( \mathbb{E}\left( z^4_{j,\ell,{\bf k}}\right) \right)^{1/2}\\
& \le
&  Q_3^{1/2}\sup_{j\in \{0,...,J\}}2^{-jd_*}\card(D_j) \\
& \le & Q_3^{1/2}.
\end{eqnarray*}
Therefore $(\mathrm{A}1)$ is satisfied.

Let's now turn to $(\mathrm{A}4)$. Again, the Jensen inequality yields
\begin{eqnarray*}
\lefteqn{\sum_{j=j_0}^{J_*}\sum_{\ell\in B_j}\sum_{{\bf K} \in \mathcal{A}_j}  \sum_{{\bf k} \in U_{j,{\bf K}}}\mathbb{E}\left( z_{j,\ell,{\bf k}}^2\indic{\left\lbrace  \left(\sum_{{\bf k} \in U_{j,{\bf K}}} z_{j,\ell,{\bf k}}^2\right)^{1/2} > (\lambda_*2^{\delta j}L^d)^{1/2}/2 \right\rbrace}\right)}& & \\
& \le
& \sum_{j=j_0}^{J_*}\sum_{\ell\in B_j}\sum_{{\bf K} \in \mathcal{A}_j}  \sum_{{\bf k} \in U_{j,{\bf K}} }\left(\mathbb{E}\left( z_{j,\ell,{\bf k}}^4\right)\right)^{1/2} \left( \mathbb{P}\left( \left(\sum_{{\bf k} \in U_{j,{\bf K}}} z_{j,\ell,{\bf k}}^2\right)^{1/2} > (\lambda_*2^{\delta j}L^d)^{1/2}/2 \right)\right)^{1/2} ~ .
\end{eqnarray*}
Using $(\mathrm{A}3)$, it comes that
\begin{eqnarray}\label{frop}
\lefteqn{ \sum_{j=j_0}^{J_*}\sum_{\ell\in B_j}\sum_{{\bf K} \in \mathcal{A}_j}  \sum_{{\bf k} \in U_{j,{\bf K}}}\left(\mathbb{E}\left( z_{j,\ell,{\bf k} }^4\right)\right)^{1/2} \left( \mathbb{P}\left( \left(\sum_{{\bf k} \in U_{j,{\bf K}}} z_{j,\ell,{\bf k}}^2\right)^{1/2} > (\lambda_*2^{\delta j}L^d)^{1/2}/2 \right)\right)^{1/2} } & & \nonumber \\
& \le
& C 2^{J_*(d_*+\delta+\upsilon)} Q_3^{1/2} \sup_{j\in \{j_0,...,J_*\}}\sup_{\ell\in B_j} \sup_{{\bf K} \in \mathcal{A}_j}  \left( \mathbb{P}\left( \left(\sum_{{\bf k} \in U_{j,{\bf K}}} z_{j,\ell,{\bf k} }^2\right)^{1/2} > (\lambda_*2^{\delta j}L^d)^{1/2}/2 \right)\right)^{1/2}\nonumber \\
& \le
& C n^{r} Q_3^{1/2} \sup_{j\in \{j_0,...,J_*\}}\sup_{\ell\in B_j}\sup_{{\bf K} \in \mathcal{A}_j}  \left( \mathbb{P}\left( \left(\sum_{{\bf k} \in U_{j,{\bf K}}} z_{j,\ell,{\bf k}}^2\right)^{1/2} > (\lambda_*2^{\delta j}L^d)^{1/2}/2 \right)\right)^{1/2}.\nonumber \\
&
&
\end{eqnarray}

To bound the probability term, we introduce the Borell inequality. For further details about this inequality, see, for instance, \citep{adler}.
\begin{lemma}[The Borell inequality] \label{cirel} Let $\mathcal{D}$ be a subset of $\mathbb{R}$ and $(\eta_t)_{t\in \mathcal{D}}$
be a centered Gaussian process.
Suppose that $$\mathbb{E}\left(\sup_{t\in \mathcal{D}}\eta_t\right)\le  N \, \ \ \ {and} \, \ \ \ \sup_{t\in \mathcal{D}} \mathbb{V} (\eta_t)\le  Z. $$
Then, for any $x>0$, we have
\begin{equation*}
\mathbb{P}\left(\sup_{t\in \mathcal{D}}\eta_t\ge x+N\right)\le  \exp(-{x^2}/{(2Z)}).
\end{equation*}
\end{lemma}
Let us consider the set $\mathcal{S}_2$ defined by $\mathcal{S}_2= \lbrace a=(a_{\bf k})\in \mathbb{R}^*; \ \sum_{{\bf k}\in U_{j,{\bf K}}}a_{\bf k}^{2}\le  1 \rbrace$,
and the centered Gaussian process $\mathcal{Z}(a)$ defined by
$$\mathcal{Z}(a) =\sum_{{\bf k} \in U_{j,{\bf K}}}a_{\bf k}z_{j,\ell,{\bf k}} .$$
We have by the Cauchy-Schwartz inequality $$\sup_{a\in \mathcal{S}_2} \mathcal{Z}(a) =\sup_{a \in\mathcal{S}_2}\sum_{{\bf k} \in U_{j,{\bf K}}}a_{\bf k}z_{j,\ell,{\bf k}}=\left(\sum_{{\bf k} \in U_{j,{\bf K}}}z^2_{j,\ell,{\bf k}}\right)^{1/2} . $$
In order to use Lemma \ref{cirel}, we have to investigate the upper bounds for
$\mathbb{E}(\sup_{a\in \mathcal{S}_2} \mathcal{Z}(a))$ and $\sup_{a\in  \mathcal{S}_2} \mathbb{V}(\mathcal{Z}(a))$.

{\bf The upper bound for $\mathbb{E}(\sup_{a\in \mathcal{S}_2} \mathcal{Z}(a) )$.} The Jensen inequality and $(\mathrm{A}3)$ imply that
\begin{eqnarray*}
\mathbb{E}\left( \sup_{a\in \mathcal{S}_2} \mathcal{Z}(a)\right) & =& \mathbb{E}\left( \left(\sum_{{\bf k} \in U_{j,{\bf K}}} z^2_{j,\ell,{\bf k}}  \right)^{1/2} \right) \le \left(\sum_{{\bf k} \in U_{j,{\bf K}}}  \mathbb{E}( z^2_{j,\ell,{\bf k}} )\right)^{1/2}\\
& \le
&  \left(\sum_{{\bf k} \in U_{j,{\bf K}}} \left( \mathbb{E}( z^4_{j,\ell,{\bf k}} )\right)^{1/2}\right)^{1/2}  \le  Q_3^{1/4}2^{\delta j/2}L^{d/2}\\
& \le
&  Q_3^{1/4}2^{\delta j/2}(\log n)^{1/2}.
\end{eqnarray*}
So, $N=Q_3^{1/4}2^{\delta j/2}(\log n)^{1/2}$.

{\bf The upper bound for $\sup_{a\in  \mathcal{S}_2} \mathbb{V}(\mathcal{Z}(a) )$.} By assumption, for any $j\in \mathbb{N}$ and ${\bf k}\in D_j$, we have $\mathbb{E}(z_{j,\ell,{\bf k}})=0$. The assumption $(\mathrm{A}4)$ yields
\begin{equation*}
\sup_{a\in \mathcal{S}_2}\mathbb{V} \left(\mathcal{Z}(a)\right)  =  \sup_{a\in \mathcal{S}_2}\mathbb{E}\left( \left( \sum_{{\bf k}\in U_{j,{\bf K}}}  a_{{\bf k}}z_{j,\ell,{\bf k}}\right)^2\right)   \le Q_4 2^{\delta j}.
\end{equation*}
It is then sufficient to take $Z=Q_4 2^{\delta j}$.

Combining the obtained expressions of $N$ and $Z$ with Lemma \ref{cirel}, for any $j\in \{j_0,...,J_*\}$, ${\bf K} \in \mathcal{A}_j$ and  ${\bf k} \in U_{j,{\bf K}}$, we have
\begin{eqnarray*}
\lefteqn{\mathbb{P}\left(\left(\sum_{{\bf k} \in U_{j,{\bf K}}}z^2_{j,\ell,{\bf k}}\right)^{1/2} \ge (\lambda_*2^{\delta j}L^d)^{1/2}/2 \right)}\nonumber \\
& =
& \mathbb{P}\left( \left(\sum_{{\bf k} \in U_{j,{\bf K}}}z^2_{j,\ell,{\bf k}}\right)^{1/2} \ge (\lambda_*^{1/2}/2-Q_3^{1/4})(2^{\delta j}L^d)^{1/2} + Q_3^{1/4}(2^{\delta j}L^d)^{1/2}  \right)\nonumber \\
& = &
\mathbb{P}\left(\sup_{a\in \mathcal{S}_2}\mathcal{Z}(a)\ge (\lambda_*^{1/2}/2-Q_3^{1/4} )(2^{\delta j}L^d)^{1/2}+N\right) \nonumber \\
& \le &
 \exp\left(-(\lambda_*^{1/2}/2-Q_3^{1/4})^2 2^{\delta j}L^d  /(2Z)\right) \le n^{-r(\lambda_*^{1/2}/2-Q_3^{1/4})^2 /(2Q_4)}.
\end{eqnarray*}
Since $\lambda_*=4\left( (2Q_4)^{1/2}+Q_3^{1/4}\right)^2$, it follows that
\begin{equation}\label{rouv2}
\mathbb{P}\left(\left(\sum_{{\bf k} \in U_{j,{\bf K}}}z^2_{j,\ell,{\bf k}}\right)^{1/2} \ge (\lambda_*2^{\delta j}L^d)^{1/2}/2 \right) \le n^{-r}.
\end{equation}
Putting \eqref{frop} and \eqref{rouv2} together, we have proved $(\mathrm{A}2)$. This ends the proof of Proposition \ref{mich}.

\section{Proof of Proposition \ref{mich2}}
The proof of this proposition is similar to the one of Theorem \ref{mini2}. The only difference is that, instead of using Lemma \ref{mal}, we use Lemma \ref{bp} below.

\begin{lemma}[\citep{caisil}]\label{bp}
Let $(v_i)_{i\in \mathbb{N}^*}$ be a sequence of real numbers, $(w_i)_{i\in\mathbb{N}^*}$ be i.i.d. $\mathcal{N}(0,1)$ and $\sigma\in \mathbb{R}^*$. Set, for any $i\in \mathbb{N}^*$,
$$u_i=v_i+\sigma w_i.$$
Then, for any $m\in\mathbb{N}^*$ and any $\gamma >1$,  the sequence of estimates $(\tilde{u}_i)_{i=1,...,m}$ defined by $\tilde{u}_i=u_i\left( 1-\gamma m\sigma^2 \left( \sum_{i=1}^m u^2_i  \right)^{-1}\right)_+$
satisfies
$$\mathbb{E}\left( \sum_{i=1}^m (\tilde{u}_i-v_i)^2\right)  \le 2 \sigma^2 \pi^{-1/2}(\gamma -1)^{-1}m^{-1/2}e^{-(m/2)(\gamma -\log \gamma -1)}+  \gamma\min \left(  \sum_{i=1}^m v_i^2 ,  \sigma^2 m\right) .$$
\end{lemma}
To clarify, if the variables $(z_{j,\ell,{\bf k}})_{j,\ell,{\bf k}}$ are i.i.d. $\mathcal{N}(0,1)$ then Lemma \ref{bp} improves the bound of the term $B_1$ appearing in the proof of Theorem \ref{mini2}.

If we analyze the proof of Theorem \ref{mini2} and we use Lemma \ref{mal} instead of Lemma \ref{bp}, we see that it is enough to determine $\lambda_*$ such that there exists a constant $Q_5>0$ satisfying
$$\sum_{j=j_0}^{J_*}\card(B_j)\card( \mathcal{A}_j) e^{-(L^q/2)(\lambda_* -\log \lambda_* -1)} \le Q_5.$$
(It corresponds to the bound of the term $B_1$ that appears in \eqref{a2}). If $\lambda_*$ is the root of $x-\log x=3$, it comes that
\begin{eqnarray*}
\sum_{j=j_0}^{J_*}\card(B_j)\card( \mathcal{A}_j) e^{-(L^q/2)(\lambda_* -\log \lambda_* -1)} &= & c_*e^{-(L^q/2)(\lambda_* -\log \lambda_* -1)}2^{J_*(d_*+\upsilon)}\\
& \le
&  C e^{-(L^q/2)(\lambda_* -\log \lambda_* -1)}n^{r} \le Q_5.
\end{eqnarray*}
Proposition \ref{mich2} is proved.

\bibliographystyle{plainnat}
\bibliography{block-curvelets}

\end{document}